\documentclass[12 pt]{article}

\usepackage{helvet} 

\usepackage{subcaption}
\usepackage{graphicx}  
\usepackage{amssymb}
\usepackage{amsmath}
\usepackage{rotating}
\newcommand{\ez}{\epsilon_0}
\usepackage{multirow}
\usepackage{soul}
\usepackage{lipsum}
              
\usepackage[margin=1.0in]{geometry}

\title{An Efficient Finite Element Solver for a Nonuniform \\ size-modified Poisson-Nernst-Planck Ion Channel Model}

\author{Dexuan Xie\thanks{Department of Mathematical Sciences, University of Wisconsin-Milwaukee, Milwaukee, WI, 53201-0413, USA
  ({\em dxie@uwm.edu}, \url{http://www.uwm.edu/\string~dxie/}). }
}

\usepackage{subcaption}

\usepackage{hyperref}
\hypersetup{
    colorlinks=true,
    linkcolor=blue,
    citecolor=blue,
    filecolor= magenta,      
    urlcolor= magenta 
}

\newcommand{\Phit}{\tilde{\Phi}}
\newcommand{\es}{\epsilon_s}
\newcommand{\ep}{\epsilon_p}
\newcommand{\emm}{\epsilon_m}
\newcommand{\rr}{{\mathbf r}}
\newcommand{\s}{{\mathbf s}}
\newcommand{\nn}{{\mathbf n}}

\newcommand{\sdiv}{{\nabla\cdot}}

\newtheorem{thm}{Theorem}[section]

\begin{document}

\maketitle

\begin{abstract}
This paper presents an efficient finite element iterative method for solving a nonuniform size-modified Poisson-Nernst-Planck ion channel (SMPNPIC) model, along with a SMPNPIC program package that works for an ion channel protein with a three-dimensional crystallographic structure and an ionic solvent with multiple ionic species. In particular, the SMPNPIC model is constructed and then reformulated by novel mathematical techniques so that each iteration of the method only involves linear boundary value problems and nonlinear algebraic systems, circumventing the numerical difficulties caused by the strong nonlinearities, strong asymmetries, and strong differential equation coupling of the SMPNPIC model. To further improve the method's efficiency, an efficient modified Newton iterative method is adapted to the numerical solution of each related nonlinear algebraic system. Numerical results for a voltage-dependent anion channel (VDAC) and a mixture solution of four ionic species demonstrate the method's convergence, the package's high performance, and the importance of considering nonuniform ion size effects. They also partially validate the SMPNPIC model by the anion selectivity property of VDAC.
\end{abstract}

%
  

\section{Introduction}

Ion channels are a class of proteins embedded in the biological membrane. They are biological devices, the ``valves''  of cells,  and the main controllers of many biological functions, including electrical signaling in the nervous system, heart and muscle contraction, and hormone secretion \cite{HilleBook2001,zheng2015handbook}. To compute macroscopic ion channel kinetics (e.g. Gibbs free energy, electric currents, transport fluxes, membrane potential, and electrochemical potential), several Poisson-Nernst-Planck (PNP) ion channel models were developed and solved numerically in the past decades (see \cite{Chen1997,Kurnikova1999,Lu:2010fk,Mathur2009,moy2000tests,ref:RouxIonChannel2004,zheng2011} for examples). We recently improved these models by introducing novel interface and boundary conditions, piecewise diffusion and bulk concentration functions, and membrane surface charge functions. To overcome the numerical difficulties caused by the solution singularities, exponential nonlinearities, multiple physical domain issues, and ionic concentration positivity constraints of a PNP ion channel model, we developed mathematical and computational techniques and used them to derive effective finite element iterative methods and corresponding software packages for solving our improved PNP ion channel models for an ion channel protein with a crystallographic three-dimensional molecular structure and an ionic solvent with multiple ionic species  \cite{XiePNPic2021,Xie4PNPicNeumann2020,Xie4PNPicPeriodic2020}. 

However, these PNP ion channel models do not consider any ion size effect in the calculation of ionic concentrations and electrostatic potential functions. Thus, they may work poorly in applications that require distinguishing ions by size. The importance of developing an ion channel model that can distinguish ions by their sizes has been well recognized in the fields of biochemistry, biophysics, and physiology. For example, if nerve cells lose their ability to distinguish sodium ions  (Na$^+$)  and potassium ions (K$^+$), which have the same charge but different sizes, life would stop in minutes. In biophysics, it is also known that an ionic solution that treats ions as spheres has quite different entropy and energy from an ionic solution that treats ions as volume-less points. To reflect ionic size effects in the calculation of electrostatic potential and ionic concentration functions, one size-modified PNP ion channel model and a related finite element solver were reported in \cite{lu2011poisson,tu2013parallel,xie2013parallel}. After more than a decade, they still represent the best work that appeared in the literature due to the challenges in modeling and computing for a molecular structure of an ion channel protein embedded in a cell membrane. Lamentably, as an extension of one early size-modified Poisson-Boltzmann model reported in \cite{borukhov1997steric}, the model treats water molecules and ions as different cubes to cause voids among ions and water molecules and a modeling redundancy problem as pointed out in \cite{SMPBEic2021}. Its iterative solver is also costly and may have a divergence problem since it involves nonlinear boundary value problems (see \cite[Eq. (33)]{xie2013parallel}) and solves each related nonlinear algebraic system by a simple fix-point iterative method (see \cite[the first equation of Eq. (31)]{xie2013parallel}). Thus, further developing size-modified PNP ion channel models and related numerical solvers remains an important research topic.

For a solution of $n$ ionic species, a size-modified PNP ion channel model is a nonlinear system that consists of $n$ size-modified Nernst-Planck boundary value problems and one Poisson boundary interface problem. Its solution gives $n$  ionic concentration functions, $\{c_i\}_{i=1}^n$, and one electrostatic potential function, $u$. Here $c_i$ denotes an ionic concentration function of species $i$. Since a Nernst-Planck equation is constructed from an electrostatic solvation free energy functional, $F$, different selections of $F$ may lead to different Nernst-Planck equations. A modification of an early functional $F$ given in \cite{borukhov1997steric} was applied for the construction of a size-modified Nernst-Planck equation in \cite[Eq. (24)]{lu2011poisson}. To fix the void and redundancy problems of this early functional $F$, we modified this $F$ in \cite[Eq. (23)]{nuSMPBE2017} and applied our improved functional to the construction of a size-modified Poisson-Boltzmann equation for a protein surrounded by an ionic solvent in \cite{nuSMPBE2017} and the construction of a size-modified Poisson-Boltzmann ion channel model in \cite{SMPBEic2021}, respectively. Using this improved functional and the modeling techniques developed in our previous work \cite{XiePNPic2021,SMPBEic2021,Xie4PNPicNeumann2020,Xie4PNPicPeriodic2020}, we construct a new nonuniform ion size modified PNP ion channel model in this work. For clarity, this new model will be denoted by SMPNPIC in the remainder of this paper.

Because of SMPNPIC's complexities, developing an effective SMPNPIC finite element solver, which includes both iterative methods and related software packages, is greatly challenging. One major difficulty comes from the strong coupling of the $n$  size-modified Nernst-Planck equations, in which, each equation depends not only on the $n$ unknown ionic concentration functions but also on the gradients of these unknown functions. To overcome this difficulty, one effective technique is to transform the $n$  ionic concentration functions $\{c_i\}_{i=1}^n$ into $n$ new functions, denoted by $\{\bar{c}_i\}_{i=1}^n$, so that the $n$ size-modified Nernst-Planck boundary value problems are transformed into the $n$ new boundary value problems such that each new problem depends only on one of $\{\bar{c}_i\}_{i=1}^n$ and does not involve any gradients of $u$ and $\{c_i\}_{i=1}^n$. Such a function transformation technique was used in the calculation of a semiconductor steady state by Slotboom \cite{slotboom1973computer}. Due to this reason, each $\bar{c}_i$ is often referred to as a Slotboom variable. This technique was applied to the SMPNP cases in \cite{lu2011poisson,xie2013parallel} and the PNP ion channel cases in \cite{XiePNPic2021,Lu:2010fk,Xie4PNPicNeumann2020,Xie4PNPicPeriodic2020}. To enable it to work for the SMPNPIC case, we start with a modification of the $n$ nonlinear algebraic equations given in \cite[Eq. (7)]{SMPBEic2021} to yield the $n$ new algebraic equations that involve $\{\bar{c}_i\}_{i=1}^n$, $\{{c}_i\}_{i=1}^n$, and $u$. From these new equations, we obtain the required function transformation and use it to derive the $n$ transformed problems. Combining these $n$ transformed problems with the $n$ new algebraic equations and a new Poisson problem, we obtain an expanded nonlinear system with $2n+1$ unknown functions --- $\{\bar{c}_i\}_{i=1}^n$, $\{{c}_i\}_{i=1}^n$, and $u$. To this end,  we can find a solution of SMPNPIC from a solution of this expanded system. 

Even so, solving this expanded system numerically is still difficult due to the singularity of $u$ at each atomic position of an ion channel protein. We overcome this singularity difficulty by adapting the solution decomposition technique developed in our prior work \cite{XiePNPic2021,SMPBEic2021,SMPBEic2019,Xie4PNPicNeumann2020,Xie4PNPicPeriodic2020}. That is, we decompose $u$ as a sum of the three potential functions, $G$, $\Psi$, and $\Phit$, which are induced from atomic charges, the potential values and charges from interfaces and boundaries, and ionic charges, respectively. Since $G$ is given in an algebraic expression and collects all the singularity points of $u$, both $\Psi$ and $\Phit$ can be found by solving linear interface boundary value problems without involving any solution singularity issues. Moreover, $\Psi$ is independent of $\Phit$ and $\{c_i\}_{i=1}^n$. Thus, we can calculate it before the search for $\Phit$ and $\{c_i\}_{i=1}^n$. In this way, we can reformulate the expanded system as a nonlinear system for computing $\Phit$, $\{c_i\}_{i=1}^n$, and $\{\bar{c}_i\}_{i=1}^n$ without involving any potential function singularity. 

We next develop a damped block iterative method for solving this nonlinear system to reduce the computational complexity and computer memory requirement.  This method stands for one major advance we made in this work. Specifically, the $2n+1$ unknown functions of the system are grouped into three blocks: Block~1 corresponds to $\{\bar{c}_i\}_{i=1}^n$, Block~2 corresponds to $\{c_i\}_{i=1}^n$, and Block~3 corresponds to $\Phit$ in their linear finite element approximations. In Block 1, we modify each transformed boundary value problem as a linear self-adjoint problem for computing one of $\{\bar{c}_i\}_{i=1}^n$ by substituting $\Phit$ and $\{c_i\}_{i=1}^n$ with their current iterative values, resulting in  $n$ independent linear problems, which can be solved quickly one by one per iteration of Block~1.  In Block~3, we derive one linear problem from a modification of the interface boundary value problem of $\Phit$ by substituting  $\{c_i\}_{i=1}^n$ with their current iterative values. Thus, each iteration of Block~3 only needs to solve one linear problem. Hence, it can also be done quickly. 

However, each iteration of Block 2 is very costly since it needs to solve a system of nonlinear algebraic equations for $\{{c}_i\}_{i=1}^n$, which we derive by substituting $\{\bar{c}_i\}_{i=1}^n$ and $\Phit$ with their current iterative values. Similarly to what was done in our previous work \cite{SMPBEic2021}, we can split this large algebraic system into $N$ independent small subsystems with $N$ being the number of interior mesh points of a finite element mesh of the solvent domain $D_s$. Since each subsystem has $n$ unknowns only, which is typically very small (such as 2, 3, or 4) in most application cases, we can quickly solve it by constructing a modified Newton iterative method. In this way, each iteration of Block~2 can be done fast to yield an updated iterative value of $\{{c}_i\}_{i=1}^n$. 

Clearly, we can solve each related linear algebraic system either directly (e.g. by the Gaussian elimination method) or iteratively (e.g. by the generalized minimal residual method using incomplete LU preconditioning (GMRES-ILU) in our current implementation). Moreover, from the classic iterative theory \cite[Section 13.5]{OR1970} it can imply that our method can converge at a fast rate since it only uses three blocks in its definition. Because our method only involves well-defined linear boundary value problems and well-defined nonlinear algebraic systems, we can use it to derive a finite element solution of SMPNPIC without involving any difficulties caused by the nonlinearity, asymmetry, singularity, coupling, and concentration positive requirement of SMPNPIC.
 
In addition, a proper selection of an initial iterate is essential to guarantee the convergence of our damped block iterative method. In \cite{SMPBEic2021}, we report a size-modified Poisson-Boltzmann ion channel (SMPBIC) model and show that a SMPBIC solution can optimally describe an equilibrium status in the sense of minimizing an electrostatic free energy functional given in \cite[Eq. (6)]{SMPBEic2021}. Interestingly, we observe that the SMPBIC model can be produced from SMPNPIC by substituting each function of $\{\bar{c}_i\}_{i=1}^n$ with the corresponding bulk concentration constants, denoted by $\{c_i^b\}_{i=1}^n$. In this sense, the SMPBIC model can be regarded as an initial approximation of SMPNPIC. Thus, we can naturally select $c_i^b$ as an initial guess to $\bar{c}_i$ and a SMPBIC solution as an initial iterate of the method. 

We implemented the damped block iterative method in Python and Fortran as a software package based on the state-of-the-art finite element library from the FEniCS project \cite{fenics-book} and the software packages that we developed to solve the PNP ion channel models \cite{XiePNPic2021,Xie4PNPicNeumann2020,Xie4PNPicPeriodic2020} and the SMPBIC models  \cite{SMPBEic2021,SMPBEic2019}. We then demonstrate the performance of this software package by doing numerical tests via a crystallographic three-dimensional molecular structure of a voltage-dependent anion channel (VDAC)  \cite{baek2017improvement} and a mixture solution of four ionic species. This VDAC protein (PBD ID: 5XD0) is in the open state conformation and has an anion selectivity property. Thus, it gives us a good test case for validating SMPNPIC. Numerical results demonstrate a fast convergence rate of our damped block iterative method, the high performance of our software package, and the importance of considering nonuniform ion size effects. They also validate SMPNPIC by the VDAC anion selectivity property. 

The rest of the paper is organized as follows.  In Section 2, we present SMPNPIC. In Section 3, we reformulate SMPNPIC into the new nonlinear system. In Section 4, we construct the damped block iterative method for solving the new system.  In Section 5, we present our iterative scheme for solving each related nonlinear algebraic system. In Section~6, we report the numerical test results. Conclusions are made in Section 7. 
 
\begin{figure}
  \center
    \includegraphics[width=0.3\textwidth]{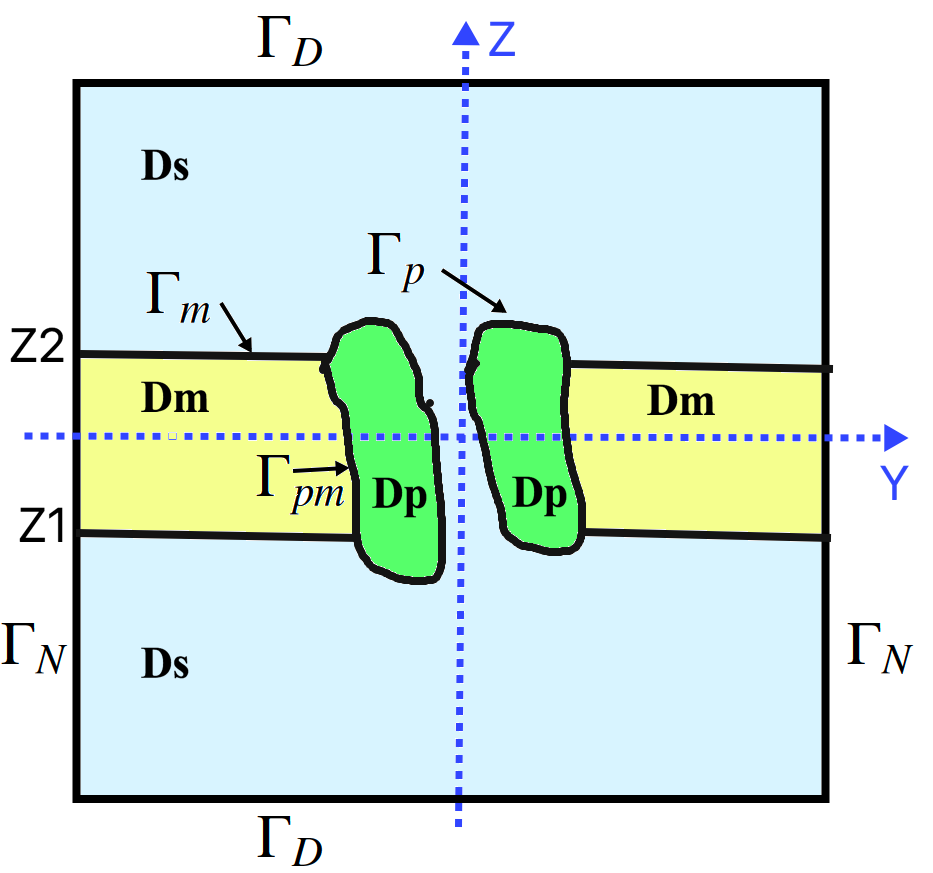}
  \caption{\small An illustration of the domain partition \eqref{DomainDecomp} and surface partition \eqref{DomainDecomp2}.}
  \label{boxdomain}
\end{figure} 

\section{A nonuniform size-modified PNP ion channel model}
In this section, we construct a nonuniform size-modified PNP ion channel model in the steady state. For clarity, we denote it by SMPNPIC. We start with a partition of an open box domain, $\Omega$,  into a protein region, $D_p$, which hosts a three-dimensional molecular structure of an ion channel protein with $n_p$ atoms, a solvent region, $D_s$, which contains a solution of $n$ ionic species, and a membrane region, $D_m$, such that
\begin{equation}\label{DomainDecomp}
\Omega = D_p\cup D_{m} \cup  D_s \cup \Gamma_p \cup \Gamma_m \cup \Gamma_{pm},
\end{equation}
where $\Gamma_p$ is an interface between $D_p$ and $D_s$, $\Gamma_m$ is an interface between $D_{m}$ and $D_s$, and $\Gamma_{pm}$ is an interface between $D_p$ and $D_m$. In particular, $\Omega$ is given by
\begin{equation}\label{Domain-omega}
  \Omega =\{ (x,y,z) \; | \: L_{x1} < x < L_{x2}, \; L_{y1} < y < L_{y2}, \; L_{z1} < z < L_{z2} \;\},
\end{equation}
 with $L_{x1}, L_{x2}, L_{y1}, L_{y2}, L_{z1}$, and $L_{z2}$ being real numbers in angstroms (\AA), the origin of the rectangular coordinate system is set at the center of the ion channel protein, the bottom and top membrane surfaces are defined by the two plane surfaces $z=Z1$ and $z=Z2$, respectively, and the normal direction of the membrane top surface coincides with the $z$-axis direction. See Figure~\ref{boxdomain} for an illustration. 

We construct SMPNPIC as a system of $n+1$ partial differential equations for computing $n$ ionic concentration functions, $\{c_i\}_{i=1}^n$, and a dimensionless electrostatic potential function, $u$. Here $c_i$ denotes the concentration function of species $i$, which is defined in the solvent region $D_s$ with the unit moles per liter (mol/L) while $u$ is defined in the box domain $\Omega$. When $u$ is given, we can get an electrostatic potential function, $\Phi$, in volts  by
\begin{equation*}\label{u2Phi}
  \Phi(\rr) =   \frac{k_{B}T}{e_{c}} u(\rr),  \quad \rr \in \Omega,
\end{equation*}
where $k_B$ is the Boltzmann constant in Joule/Kelvin, $T$ is the absolute temperature in Kelvin, and $e_{c}$ is the elementary charge in Coulomb. 

Let $\mu_i$ denote an electrochemical potential function of species $i$. It can be obtained as a variation of an electrostatic solvation free energy functional, $F$,  with respect to $c_i$. A Nernst-Planck equation in the steady state is then defined by
\begin{equation}
    \sdiv  {\cal D}_i(\rr) c_i(\rr) \nabla {\mu}_i(\rr)  =0,  \quad \rr \in D_s,  \quad i=1,2,\ldots,n,
\label{PNP2}
\end{equation}
where $\nabla$ denotes the gradient operator (i.e. $\nabla = ( \frac{\partial }{\partial x}, \frac{\partial }{\partial y}, \frac{\partial }{\partial z})$ for $\rr=(x,y,z)$) and ${\cal D}_i $ is a diffusion function of species $i$. With an electrostatic solvation free energy functional, $F$, we can derive an expression of $\mu_i$ as 
a variation of $F$ with respect to $c_i$. Thus, different $F$ may lead to different Nernst-Planck equations. Using the expression of $F$ given  \cite[Eq. (6)]{SMPBEic2021}, we can derive an expression of  $\mu_i$ as follows:
\begin{equation}
  \label{mui_eq}
        \mu_i =  k_{B}T   \gamma  \left[ Z_i u + \ln\left(\frac{c_{i}}{c_{i}^{b}}\right)   -\frac{ v_{i}}{v_{0}} \ln \big(1-   \gamma  \sum\limits_{j=1}^n v_j c_j \big) \right], \quad i=1,2,\ldots,n,
\end{equation}
where $v_i$, $Z_i$, and $c_i^b$ denote the ion size, charge number, and bulk concentration of species $i$, respectively, $v_0=\min_{1\leq i\leq n} v_i$, and $\gamma=10^{-27} N_A$ with $N_A$ denoting the  Avogadro number. Since $N_A$ is about $6.02214129\times10^{23}$ as an estimation of the number of ions per mole, we can estimate $\gamma$ by $\gamma \approx 6.022\times10^{-4}$. 
Note that all $c_i$ and the term $1 -  \gamma  \sum_{i=1}^n v_i c_i(\rr)$ are positive since the term $\gamma  \sum_{i=1}^n v_i c_i(\rr)$ gives a volume fraction of the ions in $D_s$  while the term $1 -  \gamma  \sum_{i=1}^n v_i c_i(\rr)$ gives a volume fraction of the water in $D_s$. Hence, the expression of $\mu_i$ in \eqref{mui_eq} is well defined. 

With \eqref{mui_eq}, we can find the gradient vector $\nabla \mu_i $ of $\mu_i$ in the expression 
\begin{equation*}
    \nabla \mu_i = k_{B}T   \gamma  \left[ Z_i \nabla u + \frac{1}{c_i} \nabla c_i +\frac{v_{i}}{v_0} \frac{ \gamma  \sum\limits_{j=1}^{n}v_{j} \nabla {c}_{j}(\rr)}{1 - \gamma \sum\limits_{j=1}^{n} v_j c_{j}(\rr)}  \right]. 
\end{equation*}
Substituting the above expression to \eqref{PNP2}, we derive a size-modified Nernst-Planck equation,
\begin{equation}
  \label{PNP-c_app}
\sdiv {\cal D}_i \left[ \nabla c_i(\rr) + Z_i  c_i(\rr) \nabla u(\rr)  + \frac{v_{i}}{v_0} c_i(\rr) \frac{ \gamma  \sum\limits_{j=1}^{n}v_{j} \nabla {c}_{j}(\rr)}{1 - \gamma \sum\limits_{j=1}^{n} v_j c_{j}(\rr)}  \right] =0, \quad \rr \in D_s, 
\end{equation}
and a size-modified flux vector, ${\bf J}_i$, of species $i$ in the expression
\begin{equation}
  \label{flux-def}
{\bf J}_i(\rr) = -   {\cal D}_i \left[ \nabla c_i(\rr) + Z_i  c_i(\rr) \nabla u(\rr)  + \frac{v_{i}}{v_0} c_i(\rr) \frac{ \gamma  \sum\limits_{j=1}^{n}v_{j} \nabla {c}_{j}(\rr)}{1 - \gamma \sum\limits_{j=1}^{n} v_j c_{j}(\rr)}  \right], \quad \rr \in D_s.
\end{equation}

To get boundary conditions for each size-modified Nernst-Planck equation, we split the boundary $\partial D_s$ of the solvent domain $D_s$ into three parts --- the interface part $\Gamma_p\cup \Gamma_m$, the bottom and top surface part $\Gamma_D$, and the side surface part $\Gamma_N \cap \partial D_s$. 
We then use the physical behaviors of the concentration functions on these three parts to derive three different kinds of boundary conditions as follows:

\begin{enumerate}
\item Robin boundary value conditions on the interface part $\Gamma_p\cup \Gamma_m$: 
\begin{equation}
\label{Robin-condition}
 \frac{\partial c_i(\s)}{\partial \nn_s(\s)}  +  Z_i c_i(\s)  \frac{\partial u(\s)}{\partial \nn_s(\s)} 
   + \frac{v_{i}}{v_0} c_i(\s) \frac{ \gamma  \sum\limits_{j=1}^{n}v_{j} \frac{\partial c_i(\s)}{\partial \nn_s(\s)}}{1 - \gamma \sum\limits_{j=1}^{n} v_j c_{j}(\s)} = 0,  \quad  \s\in \Gamma_p\cup \Gamma_m
\end{equation}  
for $i=1,2,\ldots, n$. Here, $\nn_s$ denotes the unit outward normal direction of the solvent region $D_s$, $\frac{\partial c_i(\s)}{\partial \nn_s(\s)} = \nabla c_i(\s) \cdot \nn_s$, and $\frac{\partial u(\s)}{\partial \nn_s(\s)} = \nabla u(\s) \cdot \nn_s$. The above equations come from the expression \eqref{flux-def} of ${\bf J}_i$ and the flux equations
\[ {\bf J}_i(\s) \cdot \nn_s(\s) = 0,  \quad \s \in \Gamma_p\cup \Gamma_m, \quad i=1,2, \ldots, n,  \]
which describes that none of the ions can penetrate the interface to enter the protein and membrane regions. 
 \item Neumann boundary value conditions on the side surface part $\Gamma_N \cap \partial D_s$:
 \begin{equation}
\label{Newmann-condition4c}
   \frac{\partial c_i(\s)}{\partial \nn_b(\s)}  = 0,  \quad \s\in \Gamma_N \cap \partial D_s, \quad i=1,2,\ldots, n,
 \end{equation}  
 where $\nn_b$ denotes the unit outward normal direction of the box domain $\Omega$. The above equations indicate that none of the ions enter the solvent region $D_s$ from the side surface of $D_s$.
 \item Dirichlet boundary value conditions on the bottom and top surface part $\Gamma_D$: 
\begin{equation}
\label{Dirichlet-condition4c}
   c_i(\s)  = g_i(\s),  \quad \s\in \Gamma_D, \quad i=1,2,\ldots,n,
 \end{equation}
 where $g_i$ is a boundary value function of species $i$. One simple selection of $g_i$ is that $g_i=c_i^b$.
 \end{enumerate}
 
 Combining each equation of \eqref{PNP-c_app} with the mixed boundary value conditions \eqref{Robin-condition}, \eqref{Newmann-condition4c}, and \eqref{Dirichlet-condition4c}, we obtain a size-modified Nernst-Planck boundary value problem as follows: 
 \begin{equation}
  \label{SMPNPIC_NP}
\left\{ \begin{array}{cl}
  \sdiv {\cal D}_i \left[ \nabla c_i(\rr) + Z_i  c_i(\rr) \nabla u(\rr)  + \frac{v_{i}}{v_0} c_i(\rr) \frac{ \gamma  \sum\limits_{j=1}^{n}v_{j} \nabla {c}_{j}(\rr)}{1 - \gamma \sum\limits_{j=1}^{n} v_j c_{j}(\rr)}  \right] =0, \quad \rr \in D_s,  &\\
\frac{\partial c_i(\s)}{\partial \nn_s(\s)}  +  Z_i c_i(\s)  \frac{\partial u(\s)}{\partial \nn_s(\s)} 
   + \frac{v_{i}}{v_0} c_i(\s) \frac{ \gamma  \sum\limits_{j=1}^{n}v_{j} \frac{\partial c_i(\s)}{\partial \nn_s(\s)}}{1 - \gamma \sum\limits_{j=1}^{n} v_j c_{j}(\s)} = 0, \quad \s\in \Gamma_p\cup \Gamma_m, &  \\
    \frac{\partial c_i(\s)}{\partial \nn_b(\s)}  = 0,  \quad \s\in \Gamma_N \cap \partial D_s, \quad
   c_i(\s)  = g_i(\s), \quad \s \in \Gamma_D, & i=1,2,\ldots, n.
\end{array}\right.
\end{equation}  

Following what was done in our previous work \cite{SMPBEic2021,Xie4PNPicNeumann2020, SMPBEic2021}, we can obtain  the following Poisson interface boundary value problem:

\begin{equation}
  \label{SMPNPIC_Poisson}
\left\{ \begin{array}{cl}
  - \epsilon_p\Delta u( \rr)  =\alpha \sum\limits_{j=1}^{n_{p}}z_{j}  \delta_{\rr_{j}}, &        \rr \in D_p,  \\
   - \epsilon_m\Delta u( \rr)  =0, &         \rr \in D_{m},  \\
 - \epsilon_s \Delta u( \rr) =   \beta \sum\limits_{i=1}^n Z_i c_i(\rr), & \rr \in D_s,  \\  
u(\s^-) = u(\s^+), \quad \ep \frac{\partial u(\s^-)}{\partial \nn_p(\s)} = \es  \frac{\partial u(\s^+)}{\partial \nn_p(\s)},   & \s\in\Gamma_p,\\
  u(\s^-) = u(\s^+), \quad  \emm  \frac{\partial u(\s^-)}{\partial \nn_m(\s)} = \es \frac{\partial u(\s^+)}{\partial \nn_m(\s)} + \tau \sigma,  & \s\in\Gamma_m,\\
  u(\s^-) = u(\s^+), \quad  \ep  \frac{\partial u(\s^-)}{\partial \nn_p(\s)} = \emm \frac{\partial u(\s^+)}{\partial \nn_p(\s)},  & \s\in\Gamma_{pm},\\
   u(\s) = g(\s),  & \s \in \Gamma_D,  \\
\frac{\partial u( \s)}{\partial \nn_b(\s)} = 0,  & \s\in \Gamma_N,
\end{array}\right.
\end{equation}  
where $z_j$ and $\rr_j$ denote the charge number and position vector of atom $j$, respectively; $\delta_{\rr_{j}}$ is the Dirac delta distribution at  $\rr_{j}$; $\sigma$ denotes a membrane surface charge density function in microcoulomb per squared centimeter ($\mu$C/cm$^2$);  $u(\s^{\pm}) = \lim_{t\rightarrow 0^+} u(\s\pm t\nn(\s))$ and 
$\frac{\partial u(\s^{\pm})}{\partial \nn(\s)} = \lim_{t\rightarrow 0^+}  \frac{\partial u(\s\pm t\nn(\s))}{\partial \nn(\s)}$ for $\nn=\nn_p, \nn_m$;
$\alpha$,  $\beta$, and $\tau$ are defined by
\begin{equation*}
\label{alpha-beta}
 \alpha = \frac{10^{10}e_{c}^{2}}{\ez k_{B}T}, \quad \beta = \frac{N_A e_{c}^{2}}{10^{17}\ez k_{B}T}, \quad \tau =  \frac{ 10^{-12} e_c}{\ez k_B T};
\end{equation*} 
and $g$ is a boundary value function. Here $\ez$ is the permittivity of the vacuum in Farad/meter and $e_{c}$ is the elementary charge in Coulomb. Using the parameter values from \cite[Table 1]{xiePBE2013}, we can estimate the three model constants $  \alpha$,  $\beta$, and $\tau$ by
\[  \alpha\approx 7042.9399, \qquad  \beta  \approx 4.2414, \qquad  \tau \approx 4.392. \]

We now {\bf obtain SMPNPIC as a nonlinear system that consists of the $n$ size-modified Nernst-Planck boundary value problems \eqref{SMPNPIC_NP} and the Poisson interface boundary value problem \eqref{SMPNPIC_Poisson}.} A solution of SMPNPIC gives the $n$ ionic concentration functions $\{c_i \}_{i=1}^n$ and electrostatic potential function $u$.

As a special case, by setting the ion sizes $v_i=0$ for $i=1,2,\ldots, n$ (i.e. without considering any ion size effect), the size-modified Nernst-Planck boundary value problems \eqref{SMPNPIC_NP} can be reduced to the following Nernst-Planck boundary value problems: 
  \begin{equation}
  \label{PNPIC_NP}
\left\{ \begin{array}{cl}
  \sdiv {\cal D}_i \left[ \nabla c_i(\rr) + Z_i  c_i(\rr) \nabla u(\rr)  \right] =0, & \rr \in D_s, \\
\frac{\partial c_i(\s)}{\partial \nn_s(\s)}  +  Z_i c_i(\s)  \frac{\partial u(\s)}{\partial \nn_s(\s)} = 0,  &  \s\in \Gamma_p\cup \Gamma_m,  \\
    \frac{\partial c_i(\s)}{\partial \nn_b(\s)}  = 0,  & \s\in \Gamma_N \cap \partial D_s,  \\
   c_i(\s)  = g_i(\s), & \s \in \Gamma_D, \quad i=1,2,\ldots,n.
\end{array}\right.
\end{equation}  
A combination of the above problems with the Poisson problem \eqref{SMPNPIC_Poisson} gives the PNP ion channel model reported in our previous work  \cite{Xie4PNPicNeumann2020}. In this sense, SMPNPIC can be regarded as an extension of this model.

From a comparison of \eqref{SMPNPIC_NP} with \eqref{PNPIC_NP}, it can be seen that each size-modified Nernst-Planck boundary value problem is not only much more strongly nonlinear and un-symmetric but also has been strongly coupled with the others. Hence, SMPNPIC is much more difficult to solve numerically than the PNP ion channel model. New mathematical and numerical techniques are required to develop an efficient iterative method for solving SMPNPIC and a related software package for computing a finite element solution of SMPNPIC on a computer. 

\section{Reformations of SMPNPIC}
In this section, we present a mathematical transformation and use it to reformulate SMPNPIC as a new nonlinear system. We then introduce a solution decomposition to avoid the singularity difficulties caused by atomic charges.

Let  $c=(c_1,c_2, \ldots, c_n)$ and $\bar{c}=(\bar{c}_1, \bar{c}_2, \ldots, \bar{c}_n)$. For a given $u$,  we construct the mathematical transformation from $c$ into $\bar{c}$ in two steps. In Step~1, we substitute the bulk concentration $c_i^b$ of the expression \eqref{mui_eq} of the size-modified electrochemical potential function $\mu_i$ with $\bar{c}_i$; we then set  $\mu_i = 0$ to derive an equation of $\bar{c}_i$ with $c$ as follows:
\begin{equation}
  \label{Slotboom-def}
   Z_i u + \ln\left(\frac{c_{i}}{\bar{c}_i}\right)   -\frac{ v_{i}}{v_{0}} \ln \big(1-   \gamma  \sum\limits_{j=1}^n v_j c_j \big)=0.
\end{equation}
In Step~2, we solve the above equation for $\bar{c}_i$ to get the mathematical transformation 
\begin{equation}
\label{out_transform}
  \bar{c}_i(\rr) =  \frac{c_{i}(\mathbf{r}) e^{ Z_{i}u(\mathbf{r})} }{ \left[ 1-  \gamma \sum\limits_{j=1}^{n}  v_{j}c_{j}(\rr)\right]^{v_i/v_0}}, \quad \rr \in D_s,  \quad i=1,2,\ldots, n.
\end{equation}
The above expression can also be written in another form of the equation of $\bar{c}_i$ with $c$:
\begin{equation}
  {c}_i(\rr) -   \bar{c}_i(\rr) \left[ 1-  \gamma \sum\limits_{j=1}^{n}  v_{j}c_{j}(\rr)\right]^{v_i/v_0} e^{ -Z_{i}u(\mathbf{r})} =0.
\label{Slotboom-eqs}
\end{equation}
Note that each $\bar{c}_i$ depends on all the ionic concentration functions due to ionic size affections. But, 
when all $v_i=0$ (i.e. ignoring size affections), the above mapping is simplified as  $n$ independent mappings from $c_i$ into $\bar{c}_i$ as follows:
\begin{equation}
\label{Slotboom-transform0}
  \bar{c}_i(\rr) = c_{i}(\mathbf{r}) e^{ Z_{i}u(\mathbf{r})}, \quad \rr \in D_s,  \quad i=1,2,\ldots, n.
\end{equation}
In this special case, $\bar{c}_i$ is often referred to as a Slotboom variable since it was introduced by Slotboom in \cite{slotboom1973computer}. Due to this reason, we can refer to each function $ \bar{c}_i$ of \eqref{out_transform} as a size-modified Slotboom variable. 

With the mathematical transformation \eqref{out_transform}, we get the following theorem.

\begin{thm}
\label{thm1}
The size-modified Nernst-Planck boundary value problem \eqref{SMPNPIC_NP} can be reformulated by the mathematical transformation \eqref{out_transform} as follows:
\begin{equation}
\label{Nernst-Planck_transformed}
\left\{ \begin{array}{cl}
 \sdiv \widehat{{\cal D}}_i(u,c) \nabla \bar{c}_i(\rr) = 0,   &  \rr\in D_s, \\
   \frac{\partial \bar{c}_i(\s)}{\partial \nn_s(\s)} = 0,  &  \s\in \Gamma_p\cup \Gamma_m,  \\
  \bar{c}_i(\s) = \bar{g}_i(\s), & \s \in \Gamma_D, \\
  \frac{\partial \bar{c}_i(\s)}{\partial \nn_b(\s)} = 0,  &  \s\in \Gamma_N \cap \partial D_s,
\end{array}\right.
\end{equation}
where  $c=(c_1,c_2,\ldots,c_n)$, $\widehat{{\cal D}}_i$ is a transformed diffusion function in the expression 
\begin{equation}
 \label{transformed_diffusion}
    \widehat{{\cal D}}_i(u,c) = {\cal D}_i(\rr)  e^{ - Z_{i} u(\rr)} \left[ 1-  \gamma \sum\limits_{i=1}^{n}  v_{i}c_{i}(\rr)\right]^{v_i/v_0}, \quad  \rr\in D_s,
\end{equation}
and $\bar{g}_i$ is a transformed boundary value function in the expression 
\begin{equation}
\label{g_bar_def}
  \bar{g}_i(\s) = \frac{g_i(\s) e^{Z_i g(\s)}}{\left[ 1-  \gamma \sum\limits_{j=1}^{n}  v_{j}g_{j}(\s) \right]^{v_i/v_0}},  \quad  \s \in \Gamma_D.
\end{equation}
\end{thm}

{\em Proof.}
Applying the gradient operation $\nabla$ to \eqref{Slotboom-def},   we can get
\[    Z_i \nabla u(\rr) + \frac{1}{c_{i}(\rr)} \nabla c_i(\rr)  - \frac{1}{\bar{c}_{i}(\rr)} \nabla \bar{c}_i(\rr) 
        + \frac{v_{i}}{v_0}  \frac{ \gamma  \sum\limits_{j=1}^{n}v_{j} \nabla {c}_{j}(\rr)}{1 - \gamma \sum\limits_{j=1}^{n} v_j c_{j}(\rr)} =0. \]
Multiplying $c_i$ on the both sides of the above equation, we can obtain
\begin{equation}
\label{c_bar_c_relationship}
  \nabla c_i(\rr) + Z_i  c_i(\rr) \nabla u(\rr)  + \frac{v_{i}}{v_0} c_i(\rr) \frac{ \gamma  \sum\limits_{j=1}^{n}v_{j} \nabla {c}_{j}(\rr)}{1 - \gamma \sum\limits_{j=1}^{n} v_j c_{j}(\rr)} 
 = \frac{c_i(\rr)}{\bar{c}_i(\rr)} \nabla \bar{c}_i(\rr).
\end{equation}
 From \eqref{out_transform} we can get the ratio $c_i / \bar{c}_i$ in the expression
 \[    \frac{c_i(\rr)}{\bar{c}_i(\rr)}  =  e^{ - Z_{i}u(\mathbf{r})} \left[ 1-  \gamma \sum\limits_{j=1}^{n}  v_{j}c_{j}(\rr)\right]^{v_i/v_0}.\]
Applying the above expression to \eqref{c_bar_c_relationship}, we get a relationship identity between $c_i$ and $\bar{c}_i$ in the equation
\begin{eqnarray}
\label{c_bar_c_relationship2}
  \nabla c_i(\rr) + Z_i  c_i(\rr) \nabla u(\rr)  + \frac{v_{i}}{v_0} c_i(\rr) \frac{ \gamma  \sum\limits_{j=1}^{n}v_{j} \nabla {c}_{j}(\rr)}{1 - \gamma \sum\limits_{j=1}^{n} v_j c_{j}(\rr)}  \nonumber \\
 =  e^{ - Z_{i}u(\mathbf{r})} \left[ 1-  \gamma \sum\limits_{j=1}^{n}  v_{j}c_{j}(\rr)\right]^{v_i/v_0} \nabla \bar{c}_i(\rr).
\end{eqnarray}
Using this identity, we can transform the first equation of \eqref{SMPNPIC_NP} into the first equation of \eqref{Nernst-Planck_transformed}. 

We next derive the boundary conditions of  \eqref{Nernst-Planck_transformed} from those of \eqref{SMPNPIC_NP}.

Since $e^{ - Z_{i}u(\mathbf{r})} [ 1-  \gamma \sum\limits_{j=1}^{n}  v_{j}c_{j}(\rr)]^{v_i/v_0} > 0$, we can use \eqref{c_bar_c_relationship2} and the Robin boundary value condition of \eqref{PNPIC_NP} on the boundary $ \Gamma_p\cup \Gamma_m$ to get the Neumann boundary value condition
\begin{equation}
\label{c_bar_Neumann}
   \frac{\partial \bar{c}_i(\s)}{\partial \nn_s(\s)} = 0,  \quad  \s\in \Gamma_p\cup \Gamma_m, 
     \quad i = 1, 2, \ldots, n.
\end{equation}

Using \eqref{c_bar_c_relationship2} and the Neumann boundary value conditions $ \frac{\partial u(\s)}{\partial \nn_b(\s)} = 0$ on $\Gamma_N$ and $\frac{\partial c_i(\s)}{\partial \nn_b(\s)}  = 0$  on $\Gamma_N \cap \partial D_s$, we also can get the Neumann boundary value conditions $\frac{\partial \bar{c}_i(\s)}{\partial \nn_b(\s)} = 0$ on $\Gamma_N \cap \partial D_s$.

Furthermore, we use \eqref{out_transform} and the Dirichlet boundary condition $u(\s)=g(\s)$ on $\Gamma_D$ to transform the Dirichlet boundary condition $c_i(\s)=g_i(\s)$ on $\Gamma_D$ of \eqref{SMPNPIC_NP} as
\begin{equation}
\label{c_bar_Dirichlet}
  \bar{c}_i(\s) = \bar{g}_i(\s),
\end{equation}
where $\bar{g}_i$ has been defined in \eqref{g_bar_def}.
This completes the proof.

\vspace{3mm}

Furthermore, we have the following theorem.

\vspace{3mm}
\begin{thm}
\label{thm2}
If $u$ and $c$ are given, then each boundary value problem of \eqref{Nernst-Planck_transformed} becomes self-adjoint and independent each other. Moreover, it has a unique solution, which is positive. 
\end{thm}

{\em Proof.}
When $u$ and $c$ are given, the transformed diffusion function $\widehat{{\cal D}}_i$ becomes a known coefficient function of the boundary value problem \eqref{Nernst-Planck_transformed}. Thus, \eqref{Nernst-Planck_transformed} becomes a linear symmetric boundary value problem with respect to $\bar{c}_i$. Hence, we can use the argument from \cite[page 27]{evans1998partial} to claim that this boundary value problem has a unique positive solution. 
This completes the proof.

\vspace{3mm}

Following our previous work \cite{Xie4PNPicNeumann2020}, we can split the electrostatic potential function $u$ into three potential functions, $G, \Psi$, and $\Phit$, such that
\begin{equation}
\label{solutionSplit}
u(\rr)=G(\rr) + \Psi(\rr) + \Phit(\rr) \quad \quad   \forall \rr\in \Omega,
\end{equation}
where $G$ is the potential induced by atomic charges, $\Psi$ is the potential induced by charges and potentials on boundaries and interfaces, and $\Phit$ is the potential induced by ionic charges. With \eqref{solutionSplit}, we can reformulate the Poisson interface boundary value problem \eqref{SMPNPIC_Poisson} to derive the algebraic expression of $G$ and the linear interface boundary value problems that define $\Psi$ and $\Phit$ as given in \cite[Eq. 13, 15, 16]{Xie4PNPicNeumann2020}. Since $G$ contains all the singularity points of $u$, these two problems do not involve any singularity points of $u$. Moreover, $\Psi$ is independent of any ionic concentrations, $\{c_i\}_{i=1}^n$. Thus, we can first calculate $\Psi$ and then treat it as a known function during a numerical calculation for $\{c_i\}_{i=1}^n$ and $\Phit$. 

Consequently, we combine \eqref{Slotboom-eqs} with \eqref{Nernst-Planck_transformed} and the linear interface boundary boundary value problem for computing $\Phit$ as given in \cite[Eq. 16]{Xie4PNPicNeumann2020} to derive an expanded nonlinear system for computing the $2n+1$ unknown functions $\Phit$, $\{c_i\}_{I=1}^n$, and $\{\bar{c}_i\}_{i=1}^n$. In this way, we have avoided the singularity difficulty caused by atomic charges. 

\section{A damped block iterative method for computing a SMPNPIC finite element solution}
In this section, we construct a damped block iterative method for solving the expanded system derived in the last paragraph of Section~3. With this method, we can obtain a finite element solution of SMPNPIC without involving any numerical difficulties caused by the strong nonlinearity, asymmetry, and singularity of SMPNPIC.

Because the box domain $\Omega$ is split into a solvent region $D_s$, a protein region, $D_p$, and a membrane region, $D_m$, the potential function $u$ is defined in $\Omega$ subject to complex interface conditions among the three subregions $D_s$, $D_p$, and $D_m$. On the other hand, each ionic concentration function, $c_i$, is defined in the solvent domain $D_s$ only, resulting in a two physical domain issue in the calculation of a SMPNPIC finite element solution. To overcome the numerical difficulties caused by the interface conditions and two physical domain issue, effective techniques were developed in our previous work \cite{XiePNPic2021,Xie4PNPicNeumann2020,Xie4PNPicPeriodic2020}. In this work, we use these techniques to generate an interface fitted irregular tetrahedral mesh, $\Omega_h$, of $\Omega$ and a tetrahedral mesh, $D_{s,h}$,  of $D_s$ such that $D_{s,h}$ is a sub mesh of $\Omega_h$. We then use these two meshes to construct two linear Lagrange finite element function spaces, $\mathcal{U}$ and $\mathcal{V}$, as two finite dimensional subspaces of the  Sobolev function spaces $H^1(\Omega)$ and $H^1(D_s)$ \cite{adams2003sobolev}, respectively. We further define the subspaces $\mathcal{U}_0$ and $\mathcal{V}_0$  of $\mathcal{U}$ and $\mathcal{V}$ by
 \[   \mathcal{U}_0 = \{ u \in \mathcal{U} \; | \:  u  =0 \mbox{ on } \Gamma_D \}, \quad 
     \mathcal{V}_0 = \{ v \in \mathcal{V} \; | \:   v  =0 \mbox{ on } \Gamma_D \}. \]
To handle the operations involving functions from both $\mathcal{U}$ and  $\mathcal{V}$, we use the restriction operator, $ {\cal R}: \mathcal{U} \rightarrow \mathcal{V}$, and prolongation operator, ${\cal P}: \mathcal{V} \rightarrow \mathcal{U}$, given in \cite{Xie4PNPicNeumann2020}. 
 
We now can obtain a finite element variational problem of \eqref{Nernst-Planck_transformed} as follows: 
Find  $\bar{c}_i \in \mathcal{V}$ satisfying $\bar{c}_i= \bar{g}_i$ on $\Gamma_D$ such that
\begin{equation}
 \label{transformed_NP}
   \int_{D_s} \widehat{{\cal D}}_i(u,c) \nabla \bar{c}_i(\rr)   \nabla v_i(\rr) d\rr = 0 \quad  \forall v_i \in  \mathcal{V}_{0}, \quad i=1,2,\ldots,n,
\end{equation}
where $\bar{g}_i$ is given in \eqref{c_bar_Dirichlet}, $\widehat{{\cal D}}_i$ is defined in \eqref{transformed_diffusion}, and $u=w+ \Phit$ with $w=G+\Psi$. We also can get a linear finite element variational problem of the problem given in \cite[Eq. 16]{Xie4PNPicNeumann2020} as follows:
{Find $\Phit \in \mathcal{U}_0$ such that } 
\begin{equation}
    \label{Phit-WeakForm}
 a(\Phit, v)  = \beta \sum\limits_{j=1}^{n}Z_{j} \int_{D_{s}}  {\cal P} c_{j}(\rr) v(\rr) d\rr    \qquad \forall v \in  \mathcal{U}_0,
\end{equation}
where $a(\cdot, \cdot)$ denotes a bilinear form defined by
 \begin{equation}
\label{WeakForm1} 
 a(\Phit, v) =  \ep\int_{D_{p}}\nabla \Phit \cdot \nabla v d\rr +  \emm\int_{D_{m}}\nabla \Phit \cdot \nabla v d\rr + \es \int_{D_{s}} \nabla \Phit \cdot \nabla v d\rr. 
\end{equation}

Let $(\bar{c}^{k}, c^{k}, \Phit^{k})$ denote the $k$th iterate of our damped three-block iterative method. Here $\bar{c}^k=( \bar{c}_1^k, \bar{c}_2^k, \ldots, \bar{c}_n^k)$ and $c^k=( c_1^k,c_2^k, \ldots, c_n^k)$, which denote the $k$th iterates in Blocks 1 and 2, respectively, while  $\Phit^{k}$ being the $k$th iterate in Block 3. We assume that $\Psi$ has been found by solving a linear finite element equation given in \cite[Eq. 29]{Xie4PNPicNeumann2020}. For clarity, we set $w=G+\Psi$. Here the algebraic expression of $G$ has been given in \cite[Eq. 13]{Xie4PNPicNeumann2020}.  For a given initial iterate $(\bar{c}^0, c^0, \Phit^0)$ with $k=0$, we calculate the $(k+1)$th iterate $(\bar{c}^{k+1}, c^{k+1}, \Phit^{k+1})$ in the following four steps: 

\begin{enumerate}
\item[] {\bf Step 1.}   Define the first block update $\bar{c}^{k+1}$ by
\begin{equation}
\label{Bcj-iterate}
    \bar{c}^{k+1}(\rr) = \bar{c}^{k}(\rr) + \omega [ \bar{p}(\rr)  - \bar{c}^{k}(\rr)], \quad \rr \in D_s, 
\end{equation}
where $\omega$ is a damping parameter with $\omega \in (0, 1)$ and $\bar{p} = (\bar{p}_1, \bar{p}_2, \ldots, \bar{p}_n)$ with $\bar{p}_i \in \mathcal{V}$ being a solution of the linear finite element  variational problem: Find  $\bar{p}_i \in \mathcal{V}$ satisfying $\bar{p}_i=\bar{g}_i$ on $\Gamma_D$ such that
\begin{equation}
\label{barC-iterate}
   \int_{D_s} \widehat{{\cal D}}_i(w + \Phit^k,c^k) \nabla \bar{p}_i \nabla v_i d\rr = 0  \quad  \forall v_i \in  \mathcal{V}_{0},
    \quad i=1,2, \ldots, n,
\end{equation}
which is a modification of \eqref{transformed_NP} by substituting $w+\Phit^k$ to $u$ and $c^k$ to $c$. Since the above $n$ problems can be solved one by one, independently, Block 1 has been decoupled into $n$ sub-blocks.
\item[] {\bf Step 2.}   Define the second block update ${c}^{k+1}$ by
\begin{equation}
\label{cj-iterate}
    {c}^{k+1}(\rr) =  {c}^{k}(\rr) + \omega \left[ {p}(\rr)  - {c}^{k}(\rr) \right], \quad \rr \in D_s,
\end{equation}
where $p=( p_1,p_2, \ldots, p_n)$ with $p_i\in \mathcal{V}$ satisfying the nonlinear algebraic system: 
\begin{equation}
 \label{EQ4p}
      p_{i}(\rr) - \bar{c}_{i}^{k+1}(\rr)  \left[1-  \gamma  \sum\limits_{j=1}^n  v_j p_j(\rr)  \right]^{\frac{ v_i}{v_0}} e^{-Z_{i} {\cal R} \left[w(\rr) + \Phit^k(\rr) \right]} =0, 
       \quad  \rr\in D_{s}, \quad  i=1,2, \ldots, n,
 \end{equation}
 which is a modification of \eqref{Slotboom-eqs} by substituting $w+\Phit^k$ to $u$ and $\bar{c}_{i}^{k+1}$ to $\bar{c}_{i}$.
 \item[] {\bf Step 3.}   Define the third block update $\Phit^{k+1}$ by
 \begin{equation}
\label{phit-iterate}
     \Phit^{k+1}(\rr) = \Phit^{k}(\rr) +  \omega \left[{q}(\rr) - \Phit^{k}(\rr) \right], \quad \rr \in \Omega,
\end{equation}
where $q$ is a solution of the linear  finite element  variational problem: 
{Find $q \in  \mathcal{U}_0$  such that} 
\begin{equation}
       \label{EQ4q}
        a(q,v)   = \beta \sum\limits_{j=1}^{n}Z_{j} \int_{D_{s}} {\cal P} {c}^{k+1}_{j}(\rr) v d\rr  \quad  \forall v  \in \mathcal{U}_0,
\end{equation}
which is a modification of \eqref{Phit-WeakForm} by substituting ${c}_{j}^{k+1}$ to ${c}_{j}$.
\item[] {\bf Step 4.}    Check the convergence: If the termination rules are satisfied:
\begin{equation}
    \label{Ite-stop}
    \|   \Phit^{k+1}  -  \Phit^{k}  \|_{\Omega} < \epsilon,  \max_{1\leq i\leq n} \|   \bar{c}^{k+1}_i -  \bar{c}^{k}_i  \|_{D_s} < \epsilon, 
       \max_{1\leq i\leq n} \|   {c}^{k+1}_i -  {c}^{k}_i  \|_{D_s} < \epsilon, 
\end{equation} 
stop the iteration and output $u=w+\Phit^{k+1}$ and $ c=c^{k+1}$ as a finite element solution of SMPNPIC; otherwise, increase $k$ by 1 and go back to Step~1. 
\end{enumerate}

\vspace{3mm}
In \eqref{Ite-stop}, we set $\epsilon=10^{-4}$ and define the function norms by
\[    \| u \|_{\Omega}= \sqrt{\int_{\Omega} u(\rr) d\rr} \quad \mbox{for } u\in \mathcal{U}, \qquad
      \| v \|_{D_s}= \sqrt{\int_{D_s} v(\rr) d\rr} \quad \mbox{for } v \in \mathcal{V}. \]
We also set the initial iterate $(\bar{c}^0, c^0, \Phit^0)$ by
\begin{equation}
\label{initial_iterates}
 \bar{c}^0 = c^b,   \quad c^0=\xi, \quad \Phit^0=q, 
 \end{equation}
where $c^b=(c^b_1, c^b_2,\ldots, c^b_n)$,  $\xi=(\xi_1, \xi_2,\ldots, \xi_n)$ with $\xi_i\in \mathcal{V}$, and $q \in  \mathcal{U}_0$ such that $(q, \xi)$ is a solution of the following nonlinear finite element variational system: 
\begin{subequations}
\label{system4nusmpbe}
\begin{eqnarray}
\label{system4nusmpbe_1}
      \xi_{i}(\rr) - c_i^b  \left[1-  \gamma  \sum\limits_{j=1}^n  v_j \xi_j(\rr)  \right]^{\frac{ v_i}{v_0}} e^{-Z_{i} {\cal R} \left[w(\rr) + q(\rr) \right]} =0, &  \rr\in D_{s}, \quad i=1, 2, \ldots, n, \\
     a(q,v)   - \beta \sum\limits_{j=1}^{n}Z_{j} \int_{D_{s}} {\cal P} \xi_{j}(\rr) v d\rr = 0  \hspace*{0.5in}  &  \forall v  \in \mathcal{U}_0.
\end{eqnarray}
\end{subequations}
From \cite{SMPBEic2021} it is known that a solution of the above nonlinear system can lead to a solution of the SMPBIC model, which describes an equilibrium status of the ion channel system. Thus, the initial iterate   \eqref{initial_iterates} is a natural selection for our damped block iterative method.  We can solve the nonlinear system \eqref{system4nusmpbe}, approximately, by using the iterative method reported in \cite{SMPBEic2021}.

With our damped block iterative method, we now can get a finite element solution of SMPNPIC by only solving the linear finite element equations \eqref{barC-iterate} and  \eqref{EQ4q} and the nonlinear algebraic system \eqref{EQ4p} without involving any numerical difficulties caused by the strong nonlinearity, asymmetry, and coupling of SMPNPIC. 

The efficiency of our iterative method can be improved sharply provided that the equations of \eqref{barC-iterate},  \eqref{EQ4q}, and \eqref{EQ4p} can be solved by efficient iterative algorithms. There exist several efficient iterative methods that we can select to solve the linear finite element equations \eqref{barC-iterate} and \eqref{EQ4q}. In this work, we select the generalized minimal residual method using incomplete LU preconditioning (GMRES-ILU) to do so by default. 

We developed an efficient modified Newton iterative scheme for solving the nonlinear algebraic system \eqref{system4nusmpbe_1} in our previous work \cite[Eq. (24)]{SMPBEic2021}. We observed that \eqref{system4nusmpbe_1} can be changed into the nonlinear algebraic system \eqref{EQ4p} by substituting the bulk concentration constant $c_i^b$ with the function $\bar{c}_{i}^{k+1}$. Thus, we can modify this scheme as an efficient scheme for solving \eqref{EQ4p}. It is this efficient scheme that has sharply improved the efficiency of our damped block iterative method. Due to its length in description, we describe this new scheme in the next section for clarity. 

\section{An iterative method for solving each related nonlinear algebraic system}
In this section, we present an iterative method for solving each nonlinear algebraic system of \eqref{EQ4p}. This scheme is a modification of the modified Newton iterative scheme given in \cite[Eq. (24)]{SMPBEic2021}. Its construction is motivated by Theorem~\ref{thm3}. 

\begin{thm}
\label{thm3}
Let $D_{s,h}$ be a tetrahedral mesh of the solvent domain $D_s$ and have the $N_h$ mesh points, $\rr^{(\mu)}$ for $ \mu =1,2,\ldots, N_h$, from $D_{s} $ and the boundary surface $ \Gamma_N \cap \partial D_s$. For a linear Lagrange finite element function space, $ \mathcal{V}$, defined on $D_{s,h}$, a solution of the nonlinear algebraic system \eqref{EQ4p} can be found equivalently from solving $N_h$ nonlinear algebraic systems as follows: For $\mu =1,2,\ldots, N_h$, 
\begin{equation}
 \label{System4meshPoints}
   p_{i,\mu} - \bar{c}_{i}^{k+1}(\rr^{(\mu)})  \left[1-  \gamma  \sum\limits_{j=1}^n  v_j p_{j,\mu}  \right]^{\frac{ v_i}{v_0}} e^{-Z_{i} {\cal R} \left[w(\rr^{(\mu)}) + \Phit^k(\rr^{(\mu)}) \right]} = 0, \quad i=1, 2, \ldots, n,
\end{equation}
where $p_{i,\mu}$ denotes a numerical approximation of the function value of $p_i(\rr)$ at $\rr=\rr^{(\mu)}$.
\end{thm}

{\em Proof.}
According to the definition of a linear Lagrange finite element function space, for $p_i \in  \mathcal{V}$, there exists a set of basis functions, $\{\phi_\mu \}_{\mu=1}^{N_h}$, such that  $p_i$ can be expressed by
\begin{equation}
 \label{p_i_formula}
         p_i(\rr) = \sum_{\mu=1}^{N_h} p_{i,\mu} \phi_{\mu}(\rr), \quad \rr \in D_s , \quad i=1,2, \ldots, n.
\end{equation}
Thus, we only need to determine the coefficient numbers  $\{p_{i,\mu}  \}_{\mu=1}^{N_h}$ in order to find the function $p_i$. We can obtain the nonlinear system of \eqref{System4meshPoints} directly from \eqref{EQ4p} by setting $\rr = \rr^\mu$ and then use them to determine the coefficient numbers   $\{p_{i,\mu}  \}_{\mu=1}^{N_h}$.  This completes the proof. \\


For clarity, we rewrite each nonlinear algebraic system of \eqref{System4meshPoints} in the vector form:
\begin{equation}
 \label{System4meshPoints2}
    \bar{F}(P_{\mu})=\mathbf{0}, \quad \mu =1,2,\ldots, N_h,
\end{equation}
where $P_{\mu} = (p_{1,\mu}, p_{2,\mu}, \ldots, p_{n,\mu} )$ and  $\bar{F} =\left(\bar{f}_{1}, \bar{f}_{2}, \ldots, \bar{f}_n \right)$ with $\bar{f}_i$ being given by
 \begin{equation}
    \label{fi-form2}
         \bar{f}_{i}(P_{\mu}) =  p_{i,\mu}- \bar{c}_{i}^{k+1}(\rr^{(\mu)}) \left[1-  \gamma  \sum\limits_{j=1}^n  v_j p_{j,\mu}  \right]^{\frac{ v_i}{v_0}} e^{-Z_{i} {\cal R} \left[w(\rr^{(\mu)}) + \Phit^k(\rr^{(\mu)}) \right]}.
 \end{equation}  
 Since each nonlinear system of \eqref{System4meshPoints2} has only $n$ unknown variables, $\{p_{i,\mu}  \}_{i=1}^{n}$, and $n$ is a small number (e.g., $n=2, 3, $ or 4) in applications, we can solve it by the Newton iterative method:
\begin{equation}
    \label{Newton4p}
    P_{\mu}^{j+1} = P_{\mu}^j  + \Upsilon_j,   \quad j=0, 1, 2, \ldots,
\end{equation} 
where $ P_{\mu}^j$ denotes the $j$-th iterate of the  Newton method, $\Upsilon_j$ satisfies Newton's equation  
\begin{equation}
    \label{Newton4r}
  J(  P_{\mu}^j)  \Upsilon_j = - \bar{F}(P_{\mu}^j),
\end{equation} 
and the initial iterate $P_{\mu}^0$ is set as $c^k$.  
Here  $J$ denotes a $n\times n$ Jacobian matrix of $\bar{F}$ with the $(i,j)$-th entry being the partial derivative $\partial \bar{f}_{i}/\partial p_{j,\mu}$ for $i,j=1,2, \ldots,n$, which are found in the expression
\begin{equation*}
\label{Jacobian-modification}
\frac{\partial \bar{f}_{i}}{\partial p_{j,\mu}} = \left\{\begin{array}{ll}
 1 + \gamma \frac{v_i^2}{v_0} \bar{c}_{i}^{k+1}(\rr^{(\mu)}) \left[1-  \gamma  \sum\limits_{j=1}^n  v_j p_{j,\mu}
       \right]^{\frac{ v_i}{v_0} - 1} e^{-Z_{i} {\cal R} \left[w(\rr^{(\mu)}) + \Phit^k(\rr^{(\mu)})\right]}, & j = i,\\
 \gamma \frac{v_i v_j}{v_0} \bar{c}_{i}^{k+1}(\rr^{(\mu)}) \left[1-  \gamma  \sum\limits_{j=1}^n  v_j p_{j,\mu}
       \right]^{\frac{ v_i}{v_0} - 1} e^{-Z_{i} {\cal R} \left[w(\rr^{(\mu)}) + \Phit^k(\rr^{(\mu)})\right]}, & j \neq i.
\end{array}\right.
\end{equation*}
To avoid a possible numerical overflow problem in the implementation, we do a modification to the above expression as follows: For an upper bound, $M$, a value of $-Z_{i} \left[w(\rr^{(\mu)}) + \Phit^k(\rr^{(\mu)})\right]$ is truncated as $M$ if the value is more than $M$. By default, we set $M=45$.
We then can solve each Newton equation of \eqref{Newton4r} directly by the Gaussian elimination method. 

 As soon as $\| \Upsilon_j \| < \epsilon$, we terminate the Newton iterative process, output $P_{\mu}^{j+1}$ as a numerical solution of \eqref{System4meshPoints2}, and construct a numerical solution of the nonlinear algebraic system \eqref{EQ4p} by  \eqref{p_i_formula}. By default, we set $\epsilon=10^{-8}$. 

 \section{SMPNPIC software package and numerical results}

We developed a software package for computing a finite element solution of SMPNPIC in Python and Fortran based on the state-of-the-art finite element library from the FEniCS project \cite{fenics-book} and our finite element program packages reported in  \cite{SMPBEic2021,SMPBEic2019,Xie4PNPicNeumann2020}. This new package includes the ion channel finite element mesh program package reported in  \cite{Chao4mesh2022} to generate all the required meshes. It also includes the Gaussian elimination method and the GMRES-ILU method as needed to solve each related linear finite element equation. By default, we select the GMRES-ILU method using the absolute and relative residual error tolerances as  $10^{-8}$. We implemented the fast iterative scheme given in Section~5 to solve each nonlinear algebraic system of \eqref{EQ4p}. 

In this software package, each diffusion function, ${\cal D}_i$, is set as a differentiable function given in \cite[Eq. (9)]{Xie4PNPicNeumann2020}:
\begin{equation}
    \label{Di_es_def}
 {\cal D}_i(\rr) = \left\{\begin{array}{cl} 
                        {\cal D}_{i,b}, & z < Z1 \mbox{ or } z>Z2, \\
                        {\cal D}_{i,c}  + ({\cal D}_{i,c}-{\cal D}_{i,b}){\cal I}_t(\rr), & Z2 - \eta \leq z \leq Z2, \\
                        {\cal D}_{i,c}, & Z1 + \eta \leq z \leq Z2 - \eta, \\
                       {\cal D}_{i,c}  + ({\cal D}_{i,c}-{\cal D}_{i,b}){\cal I}_b(\rr),& Z1 \leq z \leq Z1 + \eta, 
 \end{array}\right.                                                
\end{equation}
where  ${\cal D}_{i,b}$ denotes a bulk diffusion constant of species $i$, ${\cal D}_{i,c} $ is a parameter for controlling a diffusion-limited conduction (or current) rate, ${\cal I}_b$ and ${\cal I}_t$ are two interpolation functions given in  \cite[eq. (27)]{tu2013parallel}, and $\eta$ is a parameter for adjusting the buffering region size. 
We also used the boundary value functions $g_i$ and $g$ given in \cite[Eq. (8)]{Xie4PNPicNeumann2020}, from which we can get the transformed boundary value function $\bar{g}_i$ of $g_i$ in the piecewise expression:
\begin{equation}
    \label{barg_def} 
 \bar{g}_i(\s) = \left\{\begin{array}{rl} \displaystyle \frac{c^b_{i}}{(1-\gamma \sum_{j=1}^n v_j c_{j}^b)^{v_i/v_0}} e^{Z_i u_b}
                                   & \mbox{at  $z= L_{z1}$ (bottom surface)}, \\
                  \displaystyle       \frac{c^b_{i}}{(1-\gamma \sum_{j=1}^n v_j c_{j}^b)^{v_i/v_0}} e^{Z_i u_t} & \mbox{at  $z=L_{z2}$ (top surface)},
 \end{array}\right.                                           
\end{equation}
where  $u_b$ and $u_t$ denote the extracellular and intracellular potentials, respectively, and $c_i^b$ a bulk concentration of species $i$ for $i=1,2, \ldots, n$.

\begin{figure}[h]
        \centering
         \begin{subfigure}[b]{0.3\textwidth}
                \centering
                \includegraphics[width=\textwidth]{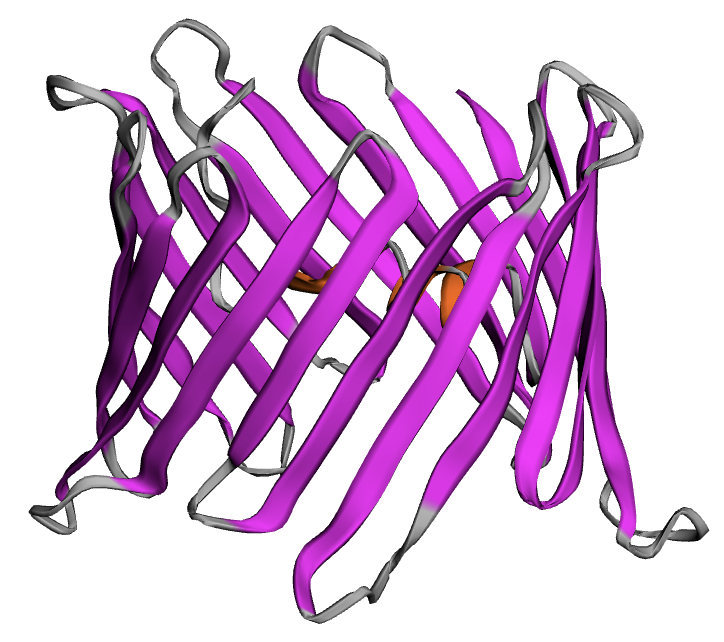}
                \vspace{1mm}
                \caption{Side view of VDAC}
        \end{subfigure}  
        \qquad \qquad \qquad
         \begin{subfigure}[b]{0.3\textwidth}
                \centering
                \includegraphics[width=\textwidth]{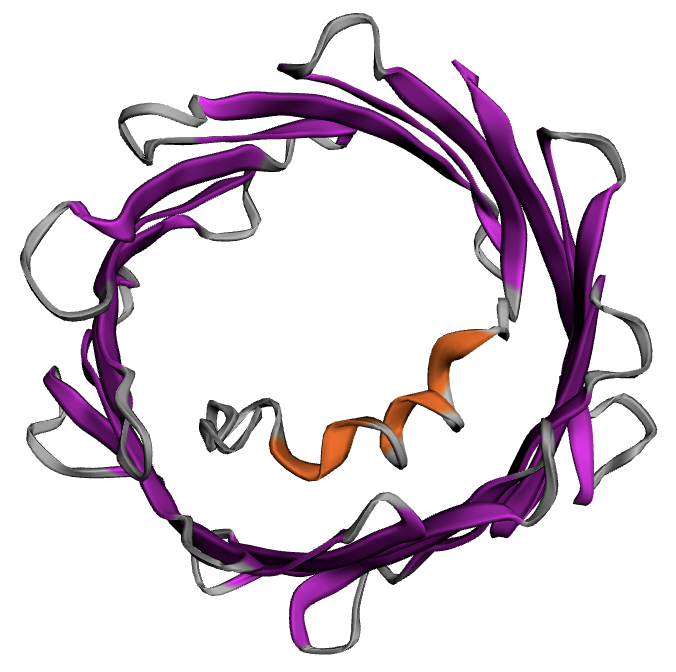}
                \caption{Top view of VDAC}
        \end{subfigure}  
              \caption{Two views of a crystallographic three-dimensional molecular structure of VDAC (PDB ID: 5XD0) in cartoon representations.
              }        
\label{structure_5XD0}  
\end{figure}

To demonstrate the performance of our SMPNPIC software package, we did numerical tests using a crystallographic molecular structure of a voltage-dependent anion channel (VDAC)  and a mixture of two salts, NaCl (sodium chloride) and KNO$_3$ (potassium nitrate), with four ionic species, Cl$^-$,   NO$_3^-$,    Na$^+$, and   K$^+$. VDAC is the most abundant protein on the outer mitochondrial membrane as the main conduit for the entry and exit of ionic species, playing a crucial role in regulating cell survival and cell death and characterizing health and diseases \cite{camara2017mitochondrial,mccommis2012role,shoshan2010vdac}. 
We got a protein data bank (PDB) file, 5XD0.pdb, of VDAC from the Orientations of Proteins in Membranes (OPM) database ({\em https://opm.phar.umich.edu}) since in OPM, each VDAC molecular structure has been transformed to satisfy our assumptions illustrated in Figure~1. From the PDB file, we got  the membrane location numbers $Z1=-11.5$ and $Z2=11.5$. We modified the PDB file as a PQR file from the PDB2PQR web server ({\em http://nbcr-222.ucsd.edu/pdb2pqr\_2.1.1/}) to get the data missed in 5XD0.pdb (such as hydrogen atoms, atomic charge numbers, and atomic radii). This PQR file is required by the mesh generation package to generate a triangular boundary surface mesh of the protein region $D_p$.  VDAC has an anion-selectivity property \cite{baek2017improvement}, which we can use to validate SMPNPIC. Two views of a cartoon backbone representation of the VDAC molecular structure are displayed in Figure~\ref{structure_5XD0}.

\begin{figure}[t]
        \centering
         \begin{subfigure}[b]{0.3\textwidth}
                \centering
                \includegraphics[width=0.9\textwidth]{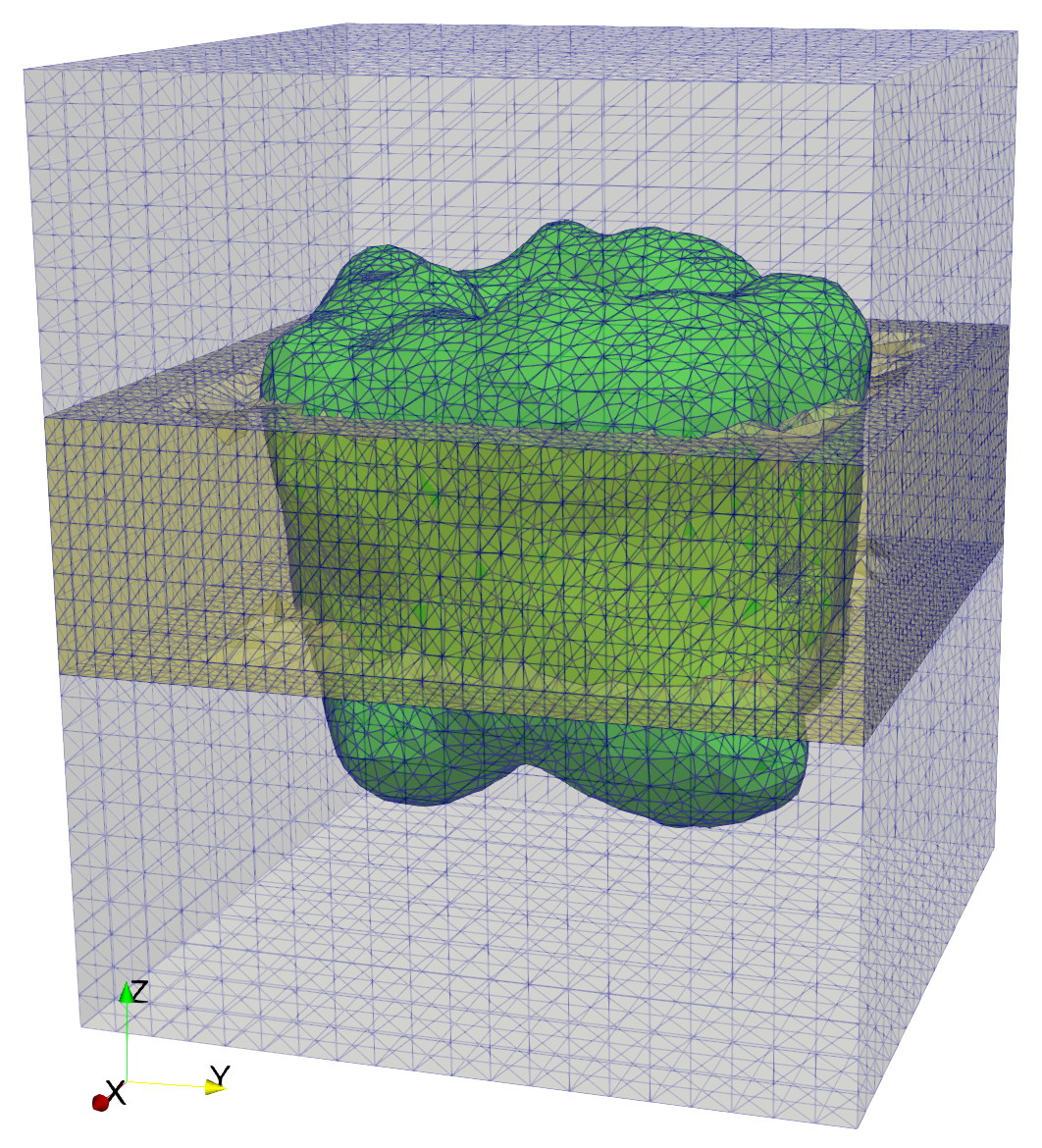}
                \caption{Side view of $\Omega_h$}
        \end{subfigure}  
        \quad 
         \begin{subfigure}[b]{0.3\textwidth}
                \centering
                \includegraphics[width=0.94\textwidth]{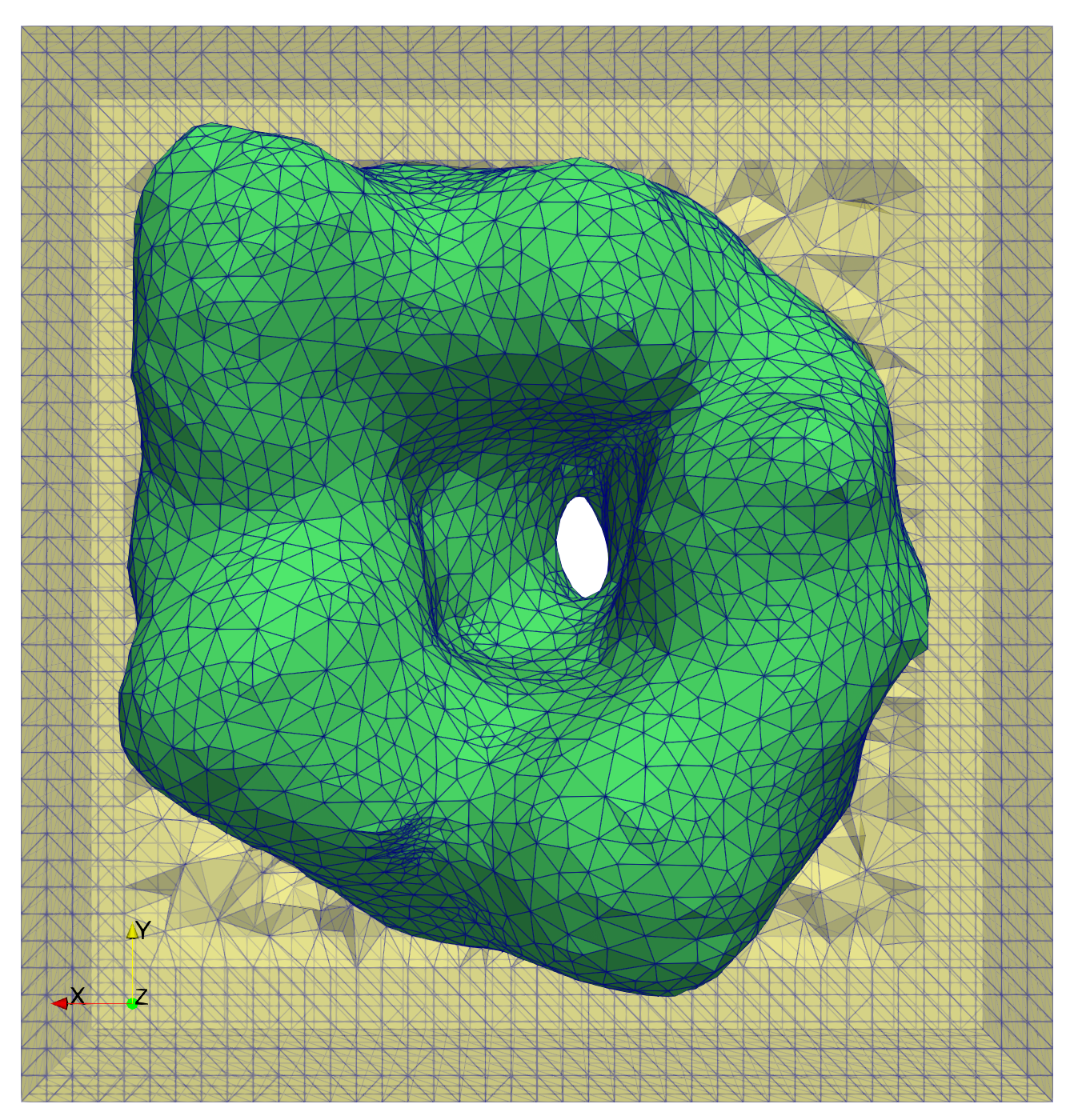}
                \caption{Top view of $\Omega_h$}
        \end{subfigure}  
        \quad 
         \begin{subfigure}[b]{0.3\textwidth}
                \centering
                \includegraphics[width=0.9\textwidth]{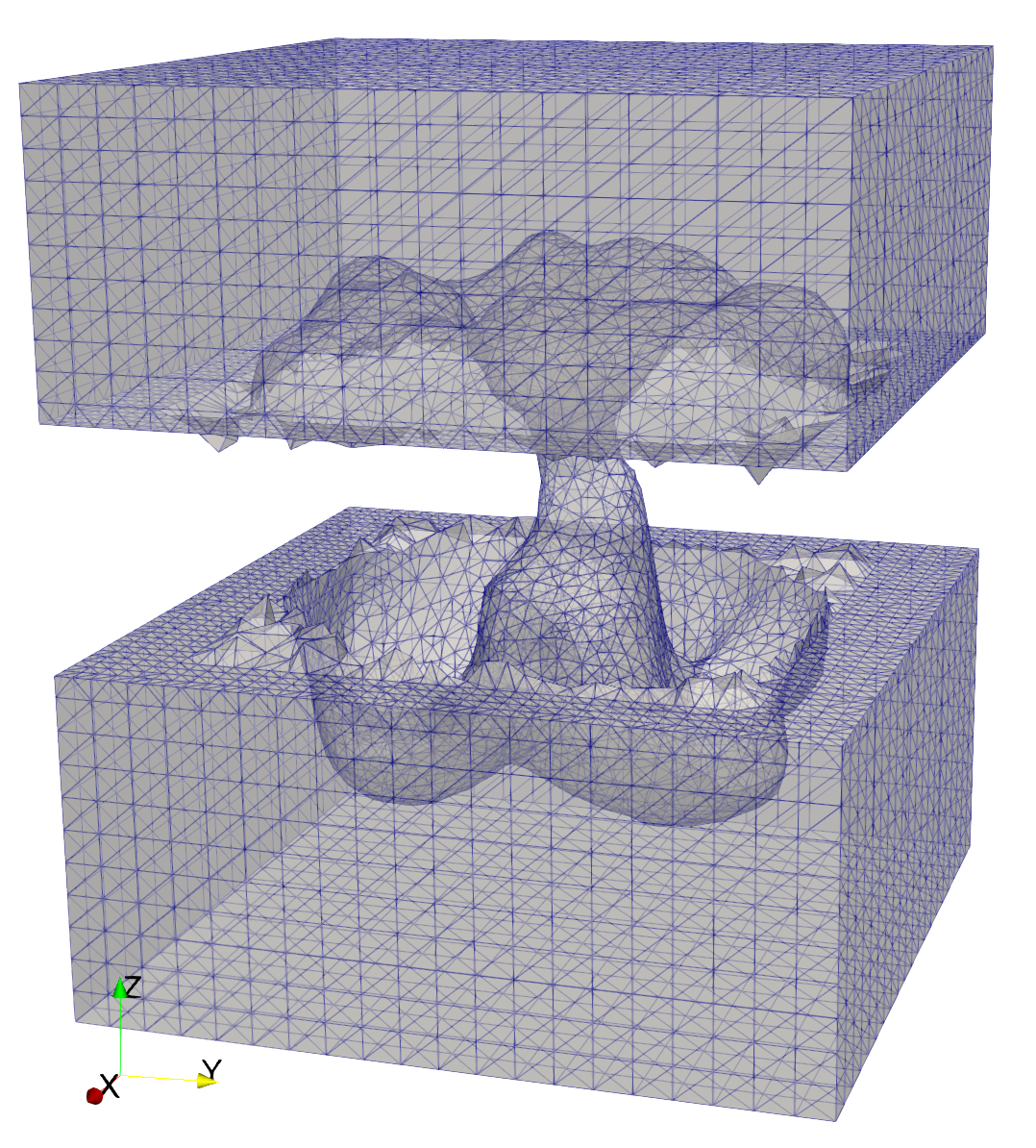}
                \caption{Side view of $D_{s,h}$ }
        \end{subfigure}
        
         \centerline{Case of Mesh 1 (coarse meshes)}
        
        \vspace{3mm}
        
         \begin{subfigure}[b]{0.3\textwidth}
                \centering
                \includegraphics[width=0.91\textwidth]{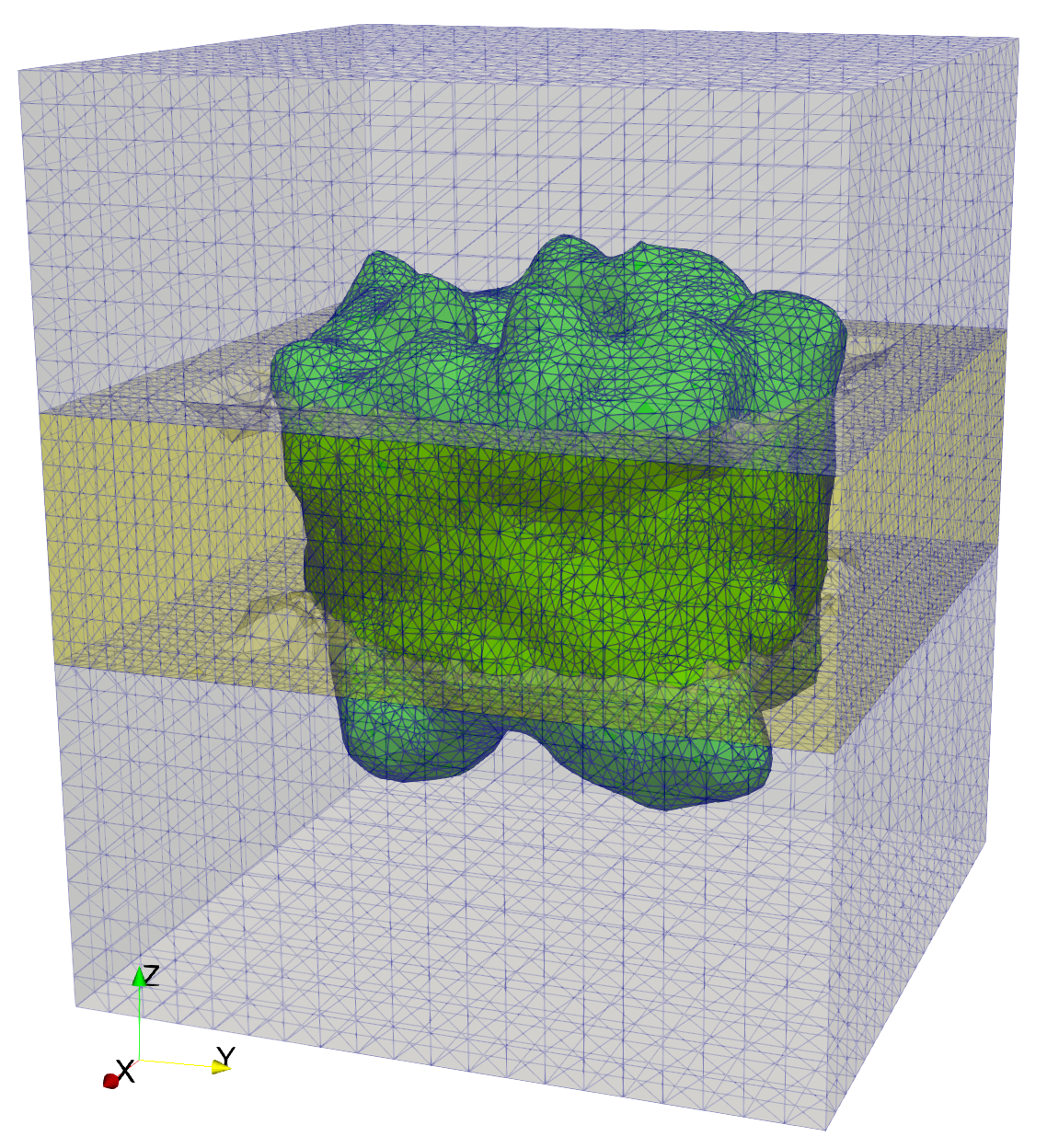}
                \caption{Side view of $\Omega_h$}
        \end{subfigure}  
        \quad 
         \begin{subfigure}[b]{0.3\textwidth}
                \centering
                \includegraphics[width=0.97\textwidth]{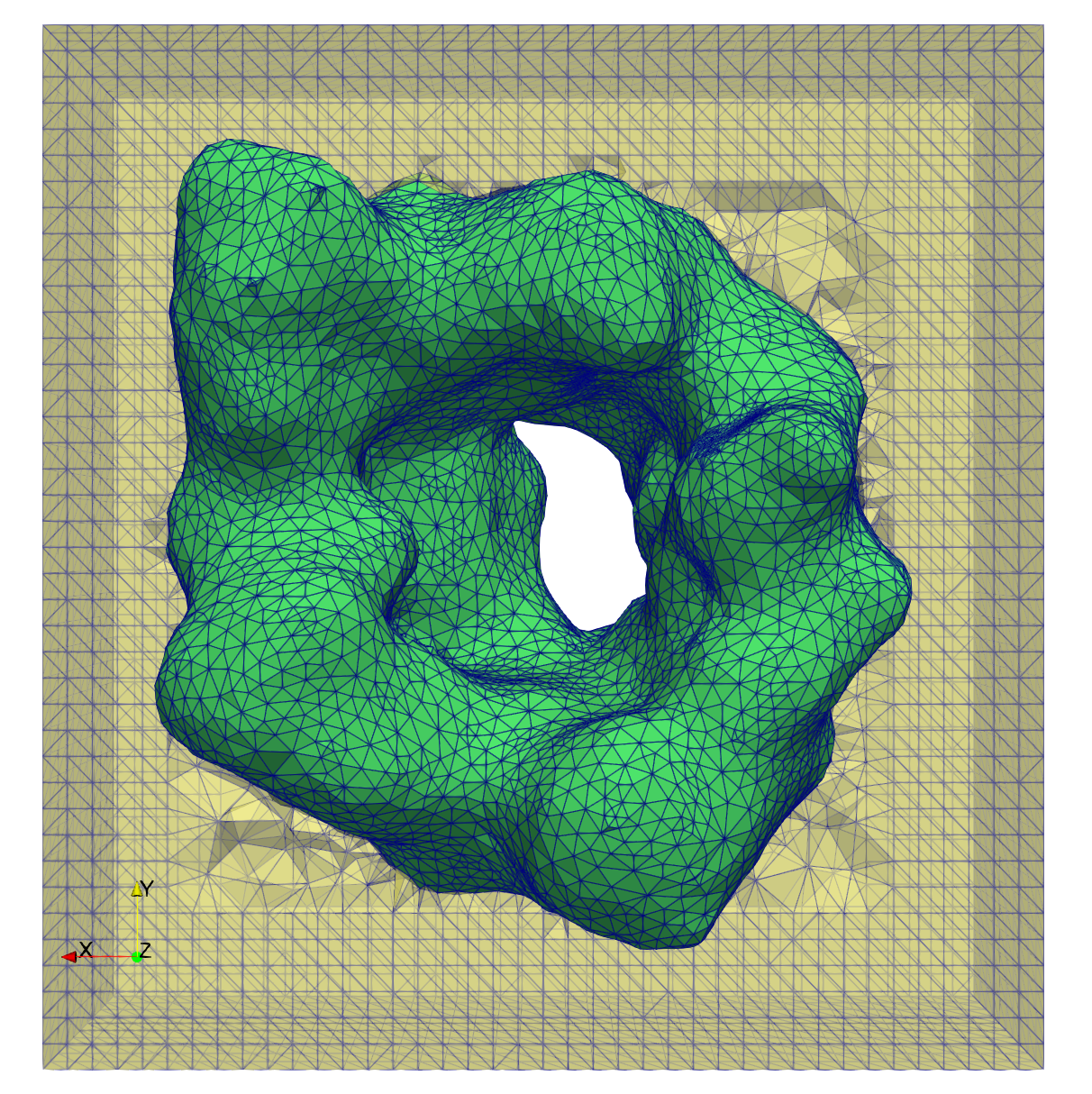}
                \caption{Top view of $\Omega_h$}
        \end{subfigure}  
        \quad 
         \begin{subfigure}[b]{0.3\textwidth}
                \centering
                \includegraphics[width=0.91\textwidth]{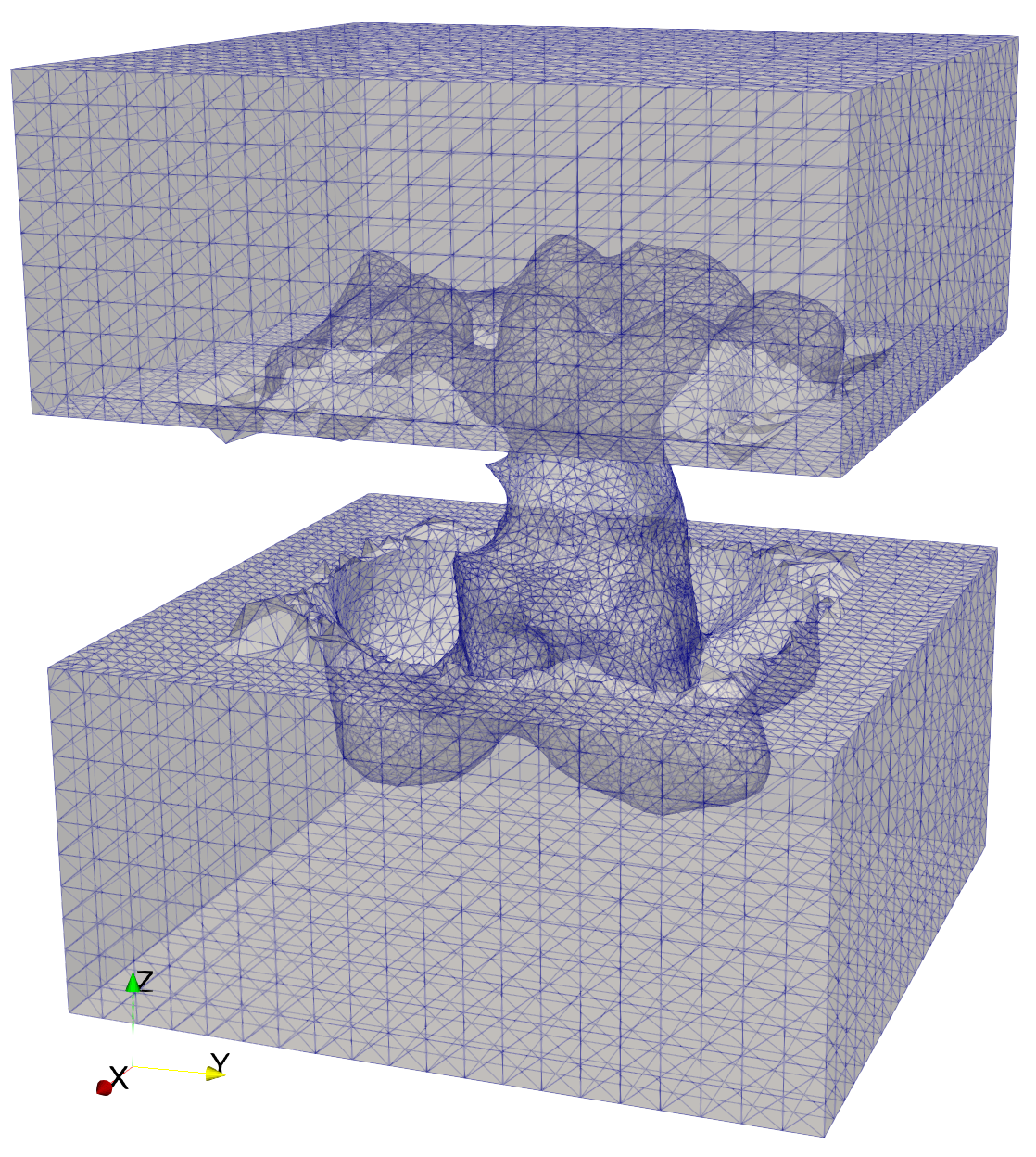}
                \caption{Side view of $D_{s,h}$ }
        \end{subfigure}
        
        \centerline{Case of Mesh 2 (fine meshes) }
        
       \caption{A comparison of the interface fitted tetrahedral box mesh $\Omega_h$ and tetrahedral solvent domain mesh $D_{s,h}$ in the case of Mesh~1 with those in the case of Mesh~2. Here the membrane region mesh $D_{m,h}$ and protein region mesh $D_{p,h}$ are colored in yellow and green, respectively, and the mesh data are given in Table~\ref{table:meshData}.}        
\label{mesh_5XD0}          
\end{figure}

We constructed a box domain, $\Omega$, of \eqref{Domain-omega} by using 
\[ L_{x1} = -59, \quad L_{x2} = 10, \quad L_{y1} = -32, \quad L_{y2} = 40, \quad  L_{z1} = -45, \quad L_{z2} = 41.\]
We then generated the interface fitted tetrahedral mesh $\Omega_h$ of $\Omega$, the meshes $D_{s,h}$, $D_{p,h}$, and  $D_{m,h}$ of the solvent region $D_{s}$, protein region $D_{p}$, and membrane region $D_{m}$ and the interface surface meshes $\Gamma_{p,h}$, $\Gamma_{m,h}$, and $\Gamma_{pm, h}$ of the interfaces $\Gamma_p$, $\Gamma_m$, and $\Gamma_{pm}$. For clarity, the set of these meshes is referred to as Mesh 1.We also generated another set of finer meshes, called Mesh 2. The major mesh data of Meshes 1 and 2 are listed in Table~\ref{table:meshData}. A comparison of Meshes~1 and 2 is given in Figure~\ref {mesh_5XD0}. 

\begin{table}[h]
\centering
\scalebox{1}{
  \begin{tabular}{|c||c|c|c|c||c|c|c|c|}
   \hline
  & \multicolumn{4}{c||}{Number of vertices  }&\multicolumn{4}{c|}{Number of tetrahedra } \\  \cline{2-9}
   &$\Omega_h$&  $D_{s,h}$ &  $D_{p,h}$ & $D_{m,h}$  &$\Omega_h$&  $D_{s,h}$ &  $D_{p,h}$ & $D_{m,h}$  \\ \hline
 Mesh 1& 26760 & 11386 & 15506 & 6323 & 155433  & 46965 & 83702 & 24766 \\   \hline
 Mesh 2& 38330   & 18672 & 21991 & 9386  & 227568  & 80312 & 109028  & 38228  \\   \hline
  \end{tabular}} 
  \caption{Mesh data  for the two mesh sets Meshes 1 and 2.}
   \label{table:meshData}
\end{table} 

Figure~\ref {mesh_5XD0} shows that the box domain mesh $\Omega_h$ has complex interfaces and the solvent domain mesh $D_{s,h}$ has an irregular complex geometrical shape. To clearly display the complex interfaces, we have highlighted the solvent, protein, and membrane meshes, $D_{s,h}, D_{p,h}$, and  $D_{m,h}$, in grey, green, and yellow colors, respectively. From Figure~\ref {mesh_5XD0}, we can see that Mesh~2 is much more irregular and complex than Mesh~1. This implies that computing a SMPNPIC finite element solution on Mesh 2 may be more difficult and more costly than that on Mesh~1. 

In the numerical tests, we used a mixed solution of 0.1 mol/L KNO$_3$ and 0.1 mol/L NaCl and ordered the four ionic species Cl$^-$,   NO$_3^-$,    Na$^+$, and   K$^+$ from 1 to 4. Thus, $c_1, c_2, c_3$, and $c_4$ denote the concentrations of Cl$^-$,   NO$_3^-$,    Na$^+$, and   K$^+$, respectively. In this case, the bulk concentration $c_i^b=0.1$ mol/L for $i=1,2,3,4$; the charge numbers  $Z_1=-1, Z_2=-1, Z_3=1$, and $Z_4=1$. We estimated the ion size $v_i$ of species $i$ by the ball volume formula, $v_i = 4\pi r_i^3/3$, where $r_i$ denotes a radius. From the website {\em https://bionumbers.hms.harvard.edu/bionumber.aspx?\&id=108517}, we got the ionic radii of Cl$^-$,   NO$_3^-$,    Na$^+$, and   K$^+$ in  \AA \; as follows:
\[ r_1=1.81, \quad r_2=2.64,\quad r_3=0.95,\quad r_4=1.33,\]
from which it gives ionic sizes in  \AA$^3 $ as follows:
\[
v_1= 24.8384, \quad v_2=77.0727, \quad  v_3=3.5914,  \quad v_4=9.8547.
\]
This ionic solution is a good selection for us to demonstrate the importance of considering nonuniform ion sizes since it contains two anions with the same charge number $-1$ (Cl$^-$ and NO$_3^-$) and two cations with the same charge number $+1$ (K$^+$ and Na$^+$) while these ions have the same balk concentration but significantly distinct sizes.

\begin{figure}[h]
        \centering
         \begin{subfigure}[b]{0.22\textwidth}
                \centering
                \includegraphics[width=\textwidth]{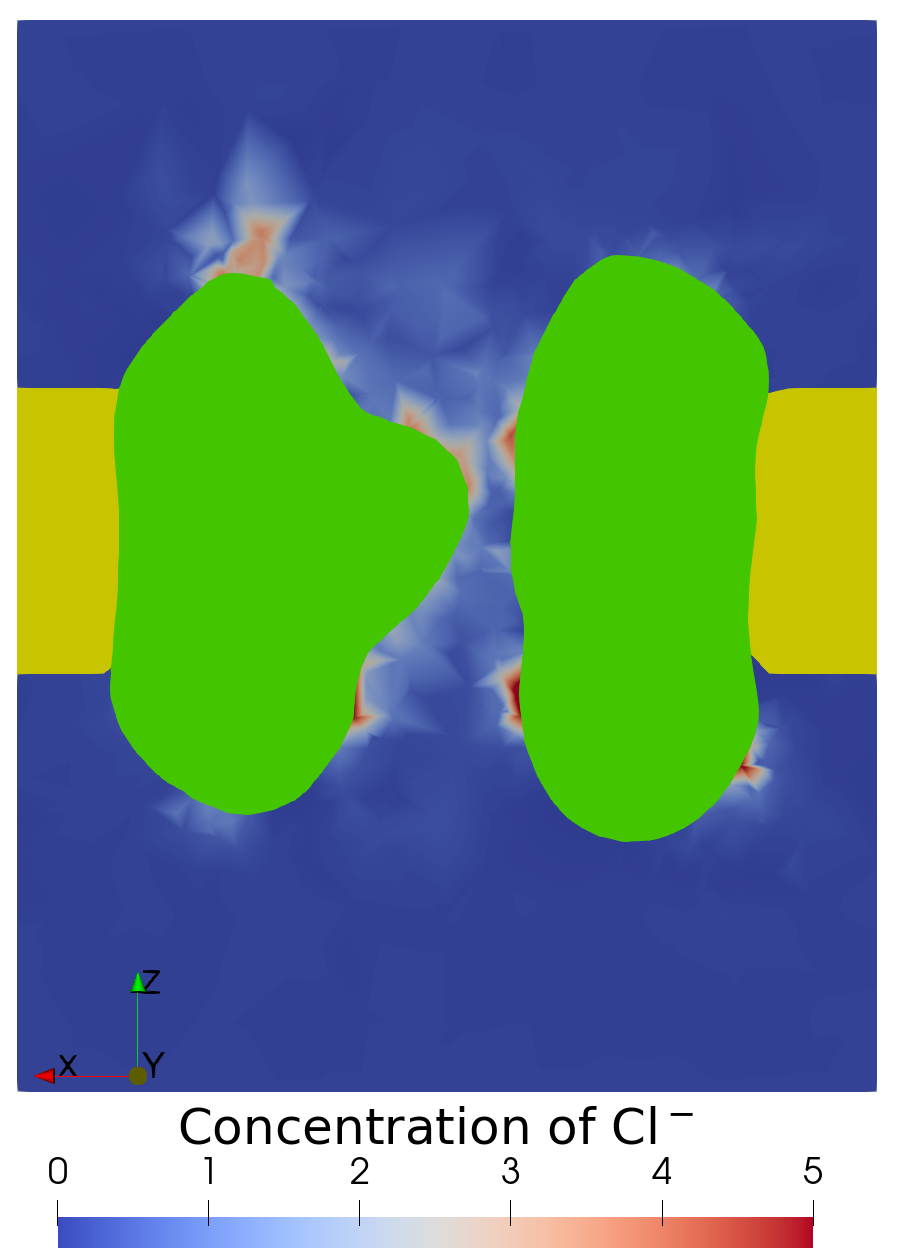}
        \end{subfigure}
        \quad
               \begin{subfigure}[b]{0.22\textwidth}
                \centering
                \includegraphics[width=\textwidth]{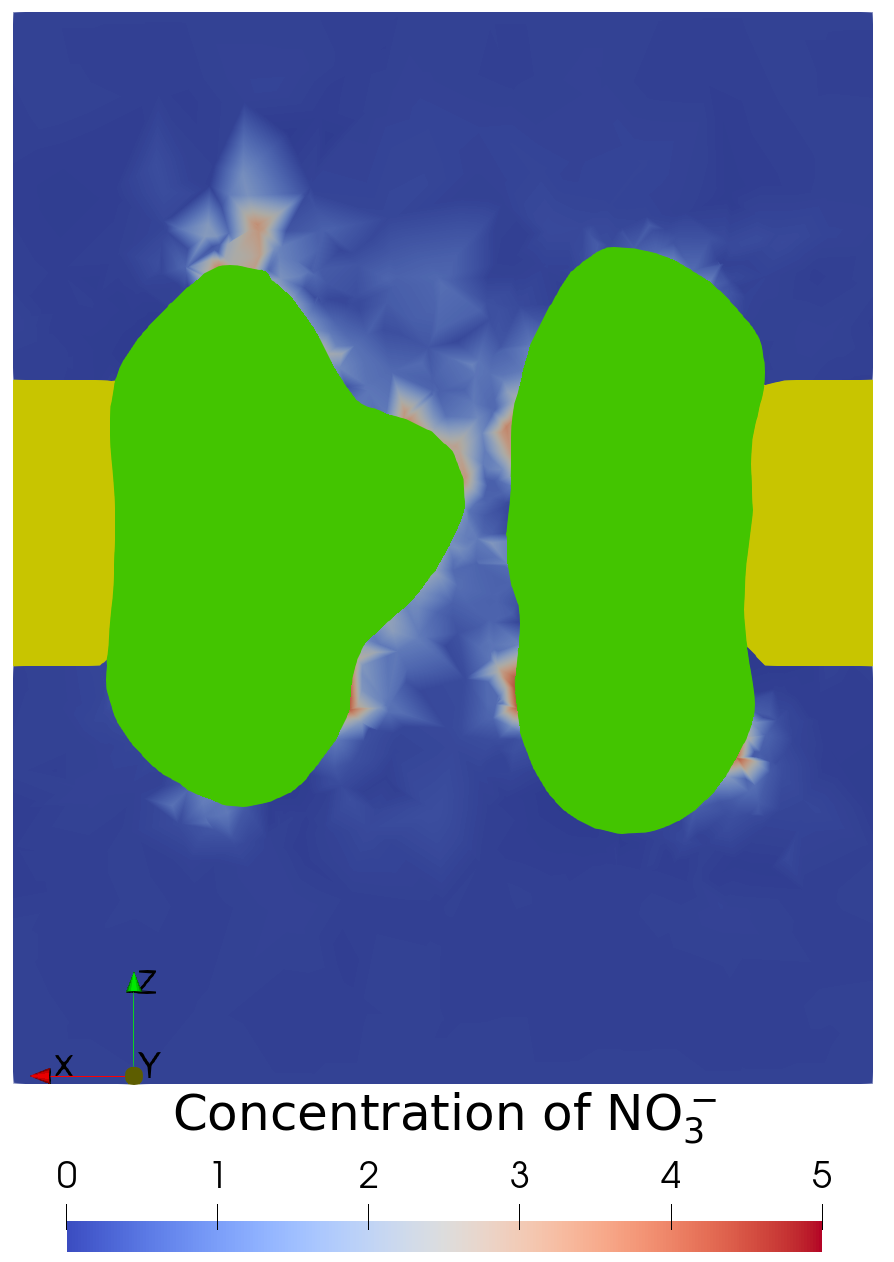}
        \end{subfigure}  
        \quad
         \begin{subfigure}[b]{0.22\textwidth}
                \centering
                \includegraphics[width=\textwidth]{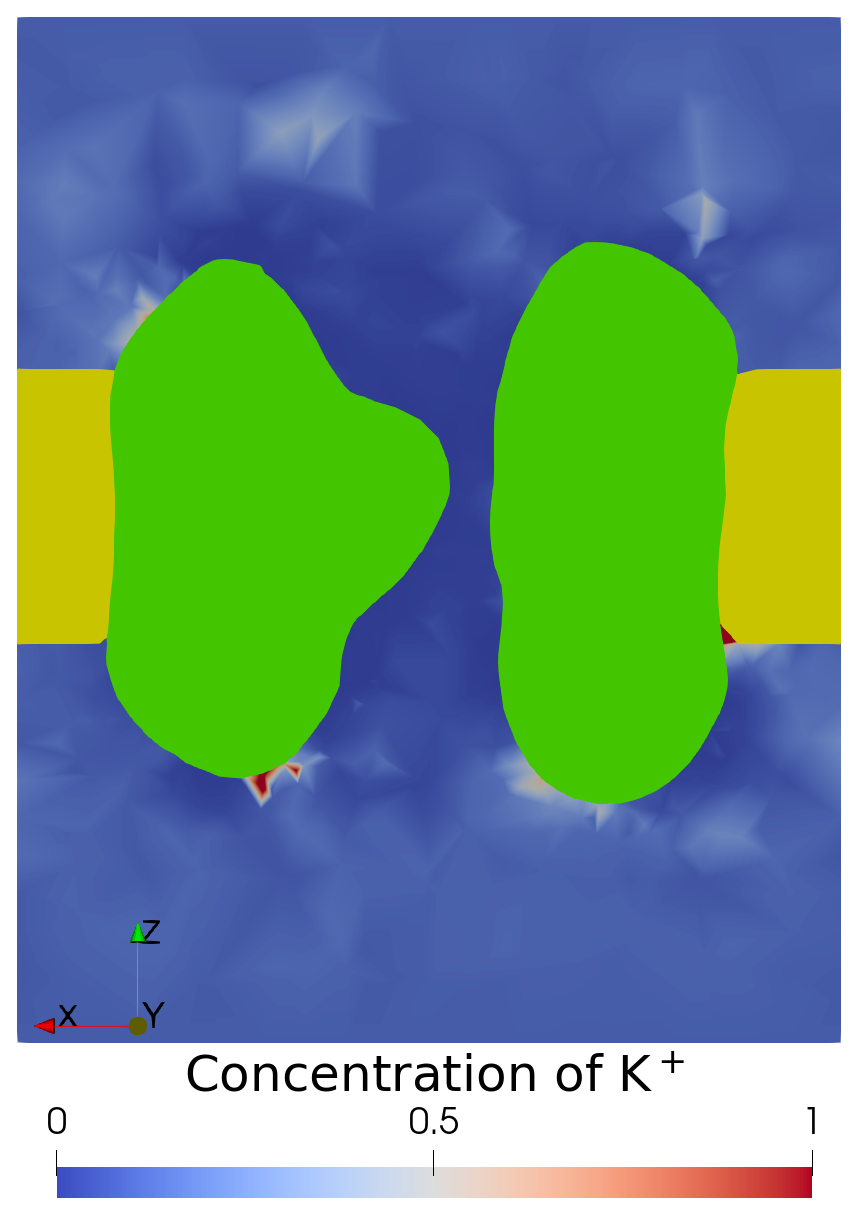}
        \end{subfigure}
        \quad
               \begin{subfigure}[b]{0.22\textwidth}
                \centering
                \includegraphics[width=\textwidth]{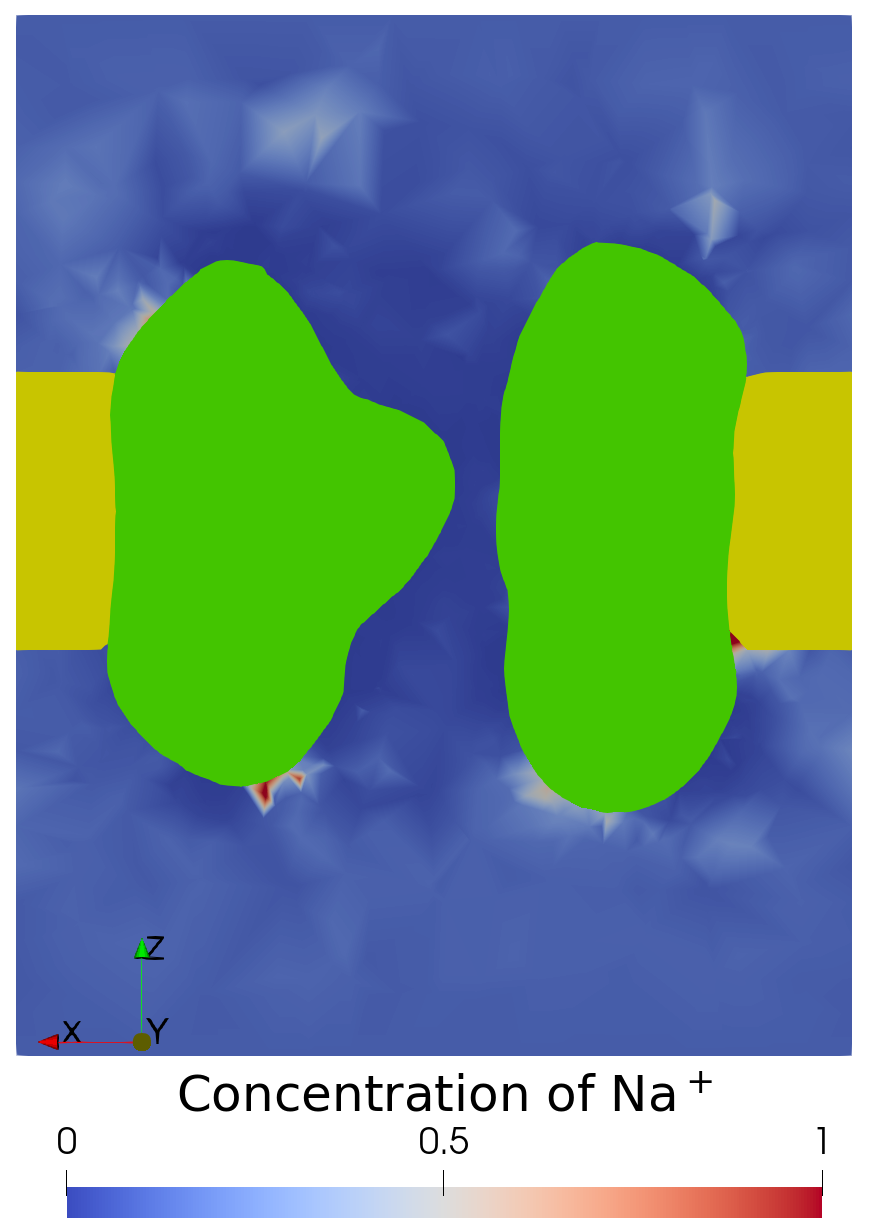}
        \end{subfigure} 
        
         {\small (a)  Case of Mesh 1}  
         
         \vspace{2mm}
         
        \begin{subfigure}[b]{0.22\textwidth}
                \centering
                \includegraphics[width=\textwidth]{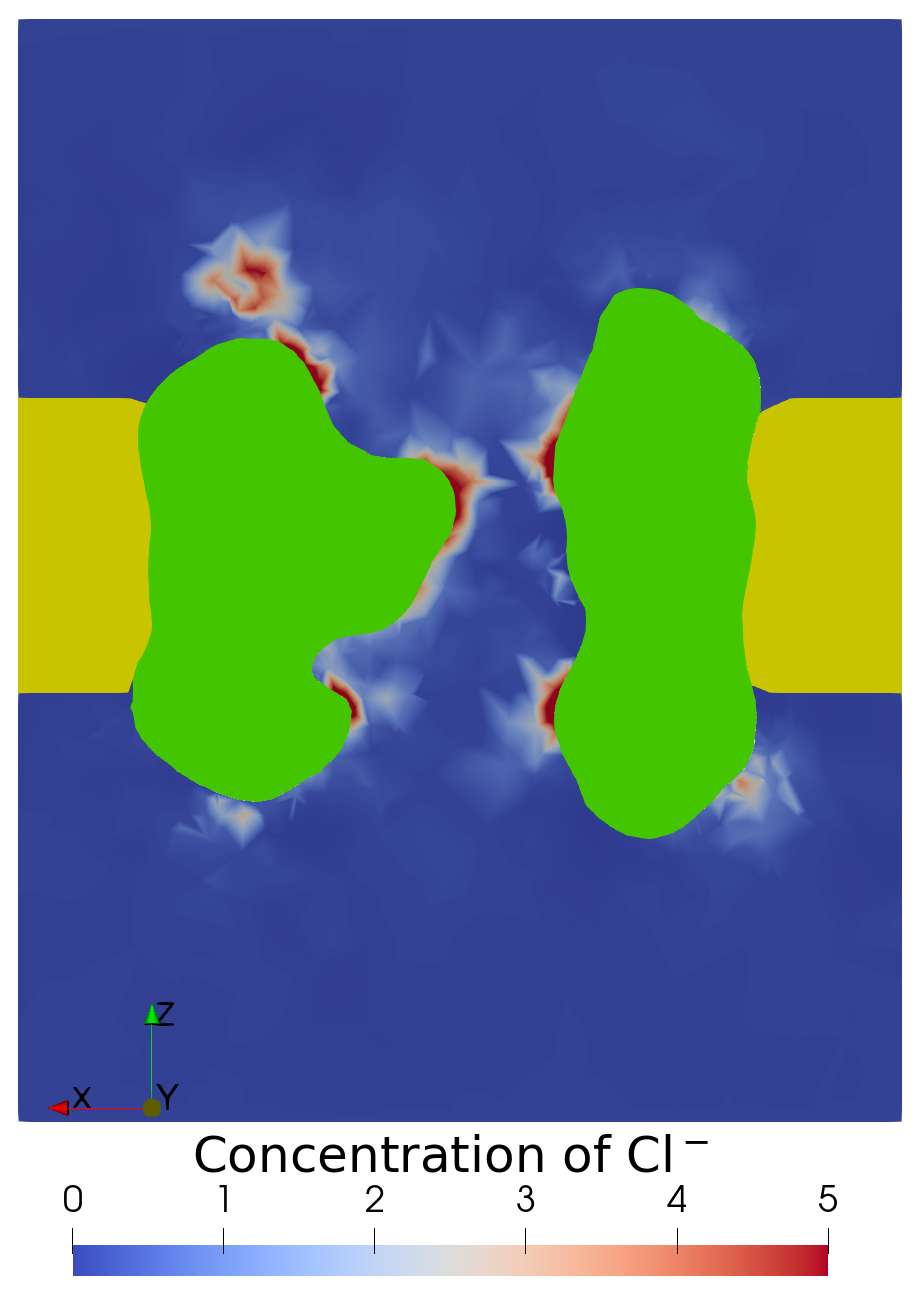}
        \end{subfigure}
        \quad
               \begin{subfigure}[b]{0.22\textwidth}
                \centering
                \includegraphics[width=\textwidth]{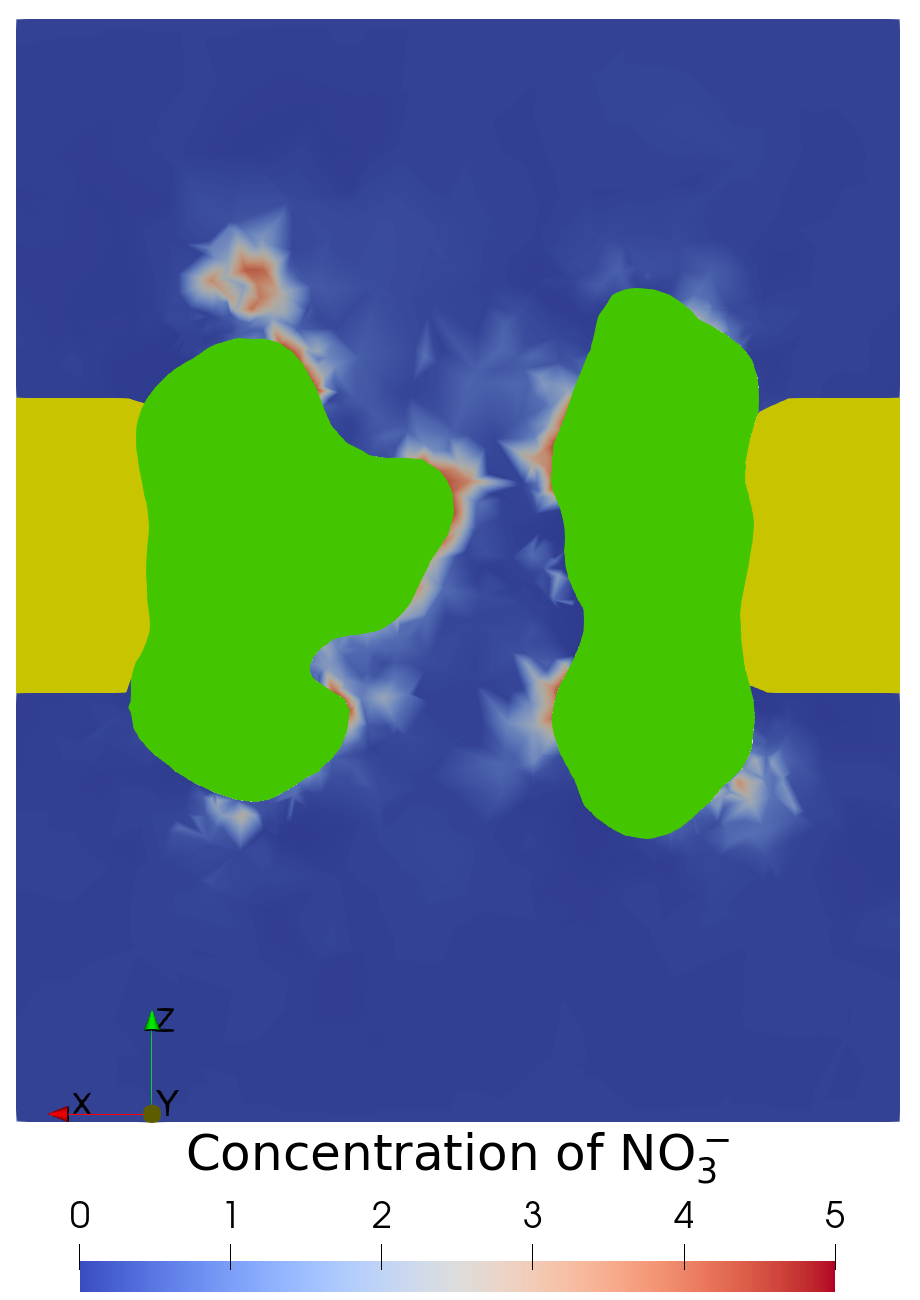}
        \end{subfigure}  
        \quad
         \begin{subfigure}[b]{0.22\textwidth}
                \centering
                \includegraphics[width=\textwidth]{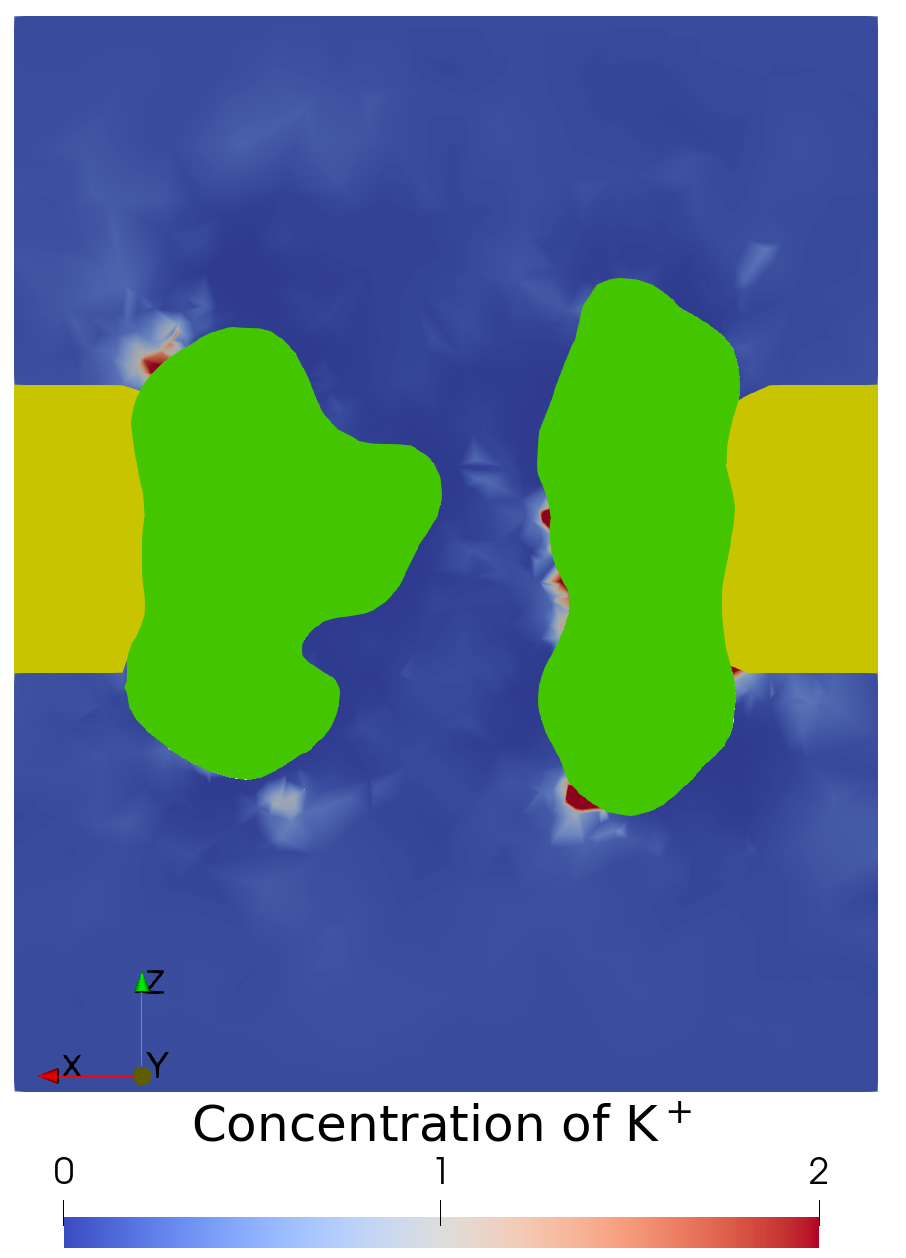}
        \end{subfigure}
        \quad
               \begin{subfigure}[b]{0.22\textwidth}
                \centering
                \includegraphics[width=\textwidth]{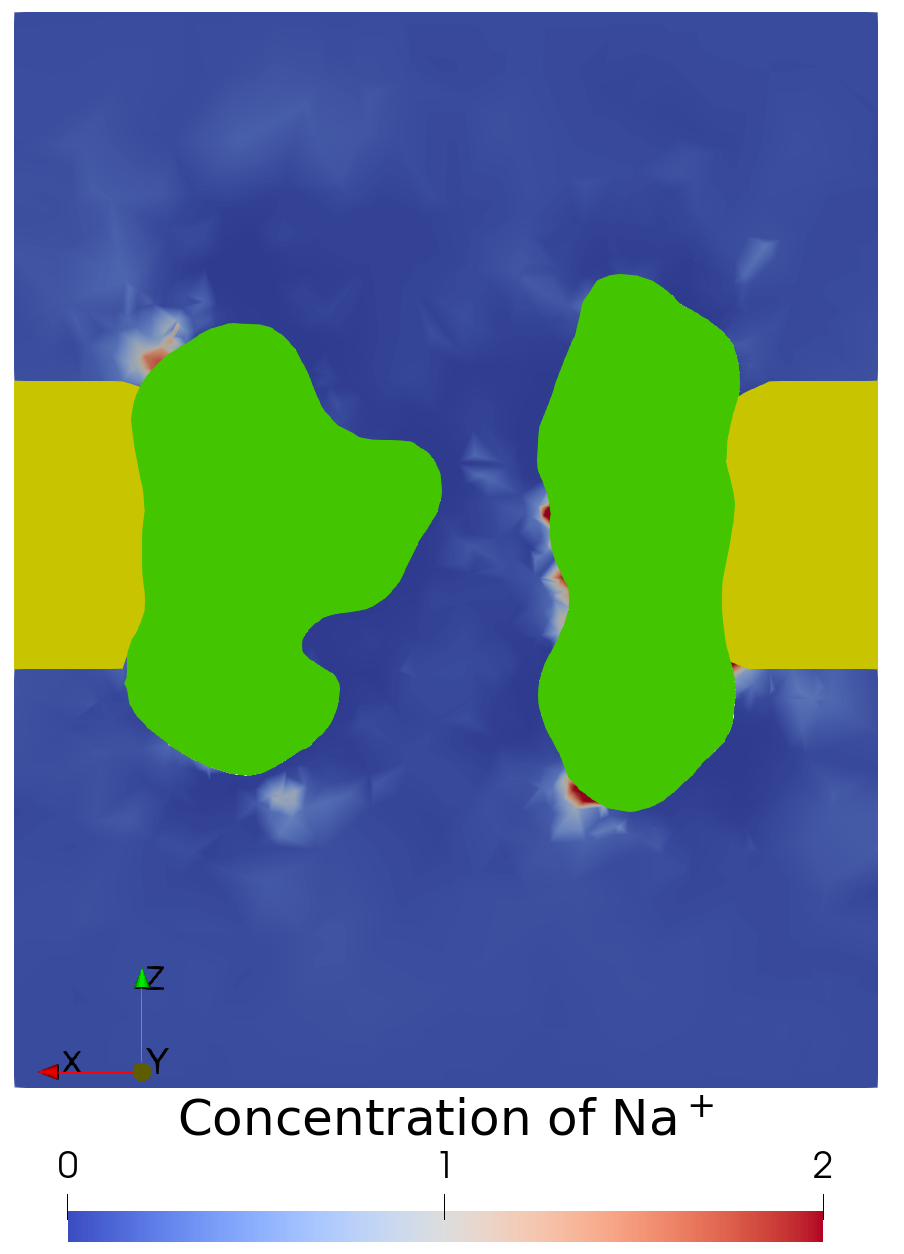}
        \end{subfigure} 
        
         {\small (b)  Case of Mesh 2}  
         
\caption{A color mapping of the ionic concentrations generated by our SMPNPIC software package on a cross-section ($y=0$) of the solvent region $D_s$ for the VDAC (PDB ID: 5XD0) and the solution of four ionic species Cl$^-$,   NO$_3^-$,    Na$^+$, and   K$^+$.}        
\label{NO3NaClK_test4}          
\end{figure}

\begin{figure}[t]
 \centering   
         \begin{subfigure}[b]{0.45\textwidth}
                \centering
                \includegraphics[width= \textwidth]{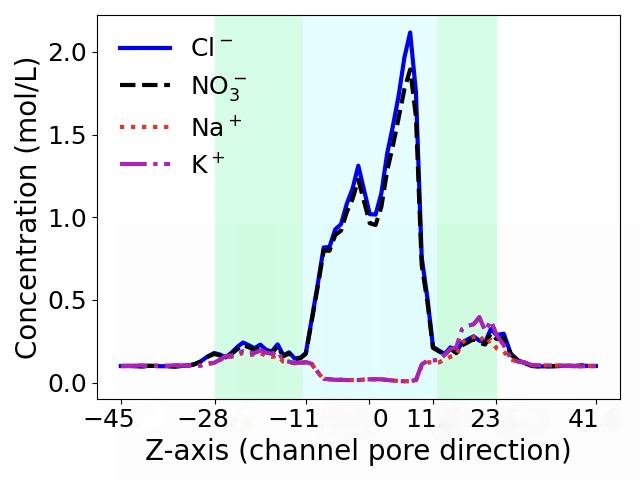}
                \caption{ Case of Mesh 1}
        \end{subfigure}
        \quad
        \begin{subfigure}[b]{0.45\textwidth}
                \centering
                \includegraphics[width= \textwidth]{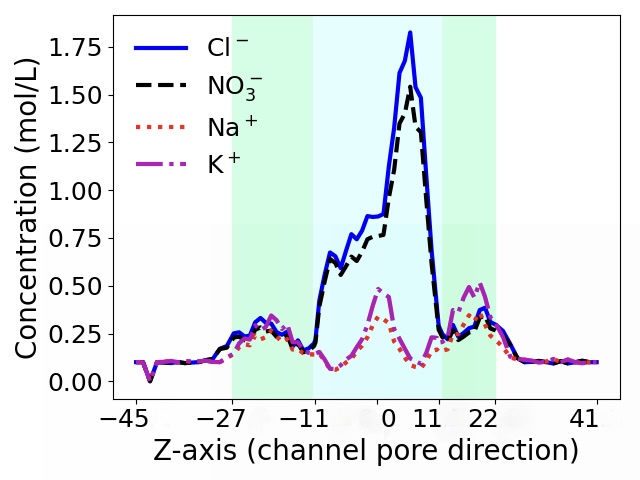}
                \caption{Case of Mesh 2}                
        \end{subfigure}     
       \caption{Mesh size influence on the ionic concentrations calculated by SMPNPIC. }         
                \label{mesh_effects}
\end{figure}

Let ${\cal D}_1^b, {\cal D}_2^b$,  ${\cal D}_3^b$, and ${\cal D}_4^b$ denote the bulk diffusion constants of  Cl$^-$,   NO$_3^-$,    Na$^+$, and   K$^+$, respectively. From the website {\em https://www.aqion.de/site/194}, we got 
\[ {\cal D}_1^b = 0.203, \quad {\cal D}_2^b = 0.190,  \quad {\cal D}_3^b = 0.133, \quad {\cal D}_4^b =0.196. \]
We then set ${\cal D}_{i,c} = \theta {\cal D}_{i,b}$ with $\theta=0.055$ to get the diffusion function ${\cal D}_i$ given in \eqref{Di_es_def}. We also set  $\ep = 2, \emm=2,$ $\es=80$, $u_t=0$, and $ u_b=0$. All the numerical tests were done on our Mac Studio computer, which has one Apple M1 Max chip and 64 GB memory. We report the test results in Tables~\ref{table:convergences_gmres} and \ref{table:performance} and Figures~\ref{NO3NaClK_test4} and \ref{mesh_effects}.

\begin{table}[h]
\centering
\scalebox{1}{
  \begin{tabular}{|c||c|c||c|c||c|c||c|c||c|c|}
   \hline
   & \multicolumn{2}{c||}{$\omega=0.30$ }&\multicolumn{2}{c||}{$\omega=0.35$  } & \multicolumn{2}{c||}{$\omega=0.38$ } & \multicolumn{2}{c||}{$\omega=0.40$ } &\multicolumn{2}{c|}{$\omega=0.41$  }   \\  \cline{2-11}
  & Ite & CPU  & Ite & CPU  & Ite & CPU & Ite & CPU & Ite & CPU  \\ \hline
\multicolumn{11}{|c|}{Solve each related linear algebraic system by  the Gaussian elimination method} \\ \cline{1-11}  
Mesh 1& 41 & 36.02& 35  & 29.14 & 32  & 35.28 & 31 & 30.31 & 30 & 26.37 \\   \hline
Mesh 2& 55   & 88.18  & 46 & 76.96 & 33 & 58.04 & 37  & 62.93 & 39 &   62.12 \\   \hline
  \multicolumn{11}{|c|}{Solve each related linear algebraic system by  the GMRES-ILU iterative method} \\ \cline{1-11}
 Mesh 1& 53 & 16.09 & 48  & 14.89  & 64  & 19.03  & 87 & 25.4268 & 77 & 21.65  \\   \hline
 Mesh 2& 109   & 57.42  & 89 & 45.63 & 104 & 52.66 & 75  & 39.01 & 74 &   38.49 \\   \hline
  \end{tabular}} 
  \caption{Convergence and performance of the damped iterative method defined in \eqref{Bcj-iterate} to \eqref{Ite-stop}. Here Ite and CPU, respectively, denote the number of iterations and computer CPU time (in seconds) determined by the iteration termination rule \eqref{Ite-stop} with $\epsilon = 10^{-4}$. }
   \label{table:convergences_gmres}
\end{table} 

Table~\ref{table:convergences_gmres} shows that the convergence and performance of our damped block iterative method, in terms of the number of iterations and the computer CPU time, were affected by the mesh size, the damping parameter $\omega$, and the iterative errors produced from the GMRES-ILU method. It also compares the convergence and performance of our damped block iterative method using the Gaussian elimination method with that using the GMRES-ILU method. Even though the iterative errors of the GMRES-ILU method disturbed the convergence rate of our damped block iterative method and caused the number of iterations to fluctuate as the damping parameter $\omega$ was increased from 0.30 to 0.41, our damped block iterative method using the GMRES-ILU method still reduced the total CPU time significantly. For example, in the test using Mesh~2 and $\omega = 0.41$, the total CPU time was reduced from about 62 seconds to about 38 seconds in spite of that the number of iterations was raised from 39 to 74. Furthermore, from Table~\ref{table:performance}, we can see that the GMRES-ILU method is more efficient than the Gaussian elimination method in the three major parts of our damped block iterative method. These test results motivated us to set the GMRES-ILU iterative method as the default linear solver in our SMPNPIC software package.

\begin{table}[h]
\centering
 \begin{tabular}{|c||c|c||c|c||c|c|}
   \hline
   &\multicolumn{2}{c||}{Calculate $\Psi$ }& \multicolumn{2}{c||}{Solve \eqref{system4nusmpbe} for  $\tilde{\Phi}^{0}$ and $c^0$  }&\multicolumn{2}{c|}{Calculate $\bar{c}$, $\tilde{\Phi}$ and ${c}$  } \\  \cline{2-7}  
 & GMRES & Direct  & GMRES & Direct & GMRES & Direct   \\ \hline
Mesh 1  & 0.28 & 2.28 &   0.19 &  2.20 &  21.65 & 26.37   \\
\hline 
 Mesh 2  & 0.42 & 3.11 & 0.27 & 3.06& 38.49  & 62.12   \\
\hline
     \end{tabular}
  \caption{The CPU time distribution (in seconds) over the three main parts of the damped iterative method. Here, GMRESS denotes the damped iterative method using the GMRES-ILU method, Direct denotes that using the Gaussian elimination method, and $\omega=0.41$. }
   \label{table:performance}
\end{table}

Figure~\ref{NO3NaClK_test4} displays the concentrations of four ionic species Cl$^-$, NO$_3^-$, K$^+$, and Na$^+$ in color mapping on the cross-section ($y=0$) of the solvent region $D_s$. Here the membrane and protein areas are colored yellow and green, respectively. From Figure~\ref{NO3NaClK_test4}, we see that the anions (Cl$^-$ and NO$_3^-$) are mostly attracted into the channel pore of VDAC while the cations (K$^+$ and Na$^+$) are expelled away from the channel pore, confirming that our SMPNPIC can retain the anion selectivity property of VDAC. From Figure~\ref{NO3NaClK_test4}, we  also see that the concentration of chloride ions (Cl$^-$) has larger values than that of nitrate ions (NO$_3^-$) within the ion channel pore area since a chloride ion is stronger in competition for space than a nitrate ion due to its smaller size (24.84 vs. 77.07 in \AA$^3$). We further see that the concentration of potassium ions (K$^+$) has larger values than that of sodium ions (Na$^+$) since a potassium ion is harder to be expelled away from the ion channel pore than a sodium ion due to its larger size (9.85 vs. 3.59 in \AA$^3$). These test results indicate that SMPNPIC can well reflect ion size effects for both anions and cations.

Figure~\ref{mesh_effects} displays the 2D mapping curves of the concentrations of  Cl$^-$,   NO$_3^-$,    Na$^+$, and   K$^+$, which we produced by using the 2D mapping scheme in \cite{Xie4PNPicNeumann2020}.
Such a 2D curve visualizes the distribution profile of an ionic species within the solvent domain $D_s$. Here the $z$-axis direction coincides with the channel pore direction; the channel pore area is between $z=-28$ and $z=23$; the membrane area is between $z=Z1$ and $z=Z2$ with $Z1=-11.5$ and $Z2=11.5$; the central part of the channel pore that links to the membrane and the other parts of the channel pore are highlighted in light-cyan and green, respectively. From Figure~\ref{mesh_effects}, we can clearly show that the concentrations of Cl$^-$ and K$^+$ have larger values than those of NO$_3^-$ and Na$^+$, especially within the central part of the channel pore, due to ionic size effects. These test results confirm that ionic sizes have significant impacts on the values of ionic concentrations. They also indicate that the anion selectivity happens mostly within the central part of an ion channel pore.

\section{Conclusions}
In this work, we have introduced a nonuniform ion size-modified Poisson-Nernst-Planck ion channel model, referred to as  SMPNPIC, and developed an efficient SMPNPIC finite element iterative method and a corresponding software package for an ion channel protein with a three-dimensional crystallographic structure and a solution with multiple ionic species. The numerical tests conducted on a voltage-dependent anion-channel (VDAC) and a solution with four ionic species demonstrate the fast convergence of the method and the high performance of the package. Moreover, these tests demonstrate the importance of considering nonuniform ion size effects and validate SMPNPIC by the anion selectivity property of VDAC. 

SMPNPIC is a system of ion size-modified Nernst-Planck equations and Poisson equations for computing multiple ionic concentration functions and one electrostatic potential function. However, due to these equations' strong nonlinearity, asymmetry, and coupling, solving SMPNPIC numerically is much more challenging than any current Poisson-Nernst-Planck ion channel models. In this work, we have overcome these difficulties through developing effective mathematical and computational techniques while overcoming the other numerical difficulties caused by solution singularities, exponential nonlinearities, multiple physical domains, and ionic concentration positivity requirements by adapting the techniques developed in our previous work. As a result, we can find a SMPNPIC finite element solution by only solving linear boundary value problems and nonlinear algebraic equation systems. Moreover, we have developed efficient iterative methods for solving these linear problems and nonlinear algebraic systems, sharply improving the efficiency of our SMPNPIC finite element solver. 

In the future, we will further improve the efficiency and usage of our SMPNPIC software package and extend its application to the calculation of various ion channel kinetics, such as Gibbs free energy change, membrane potential, and electric currents. These developments will turn our package into a highly valuable tool for biochemists, biophysicists, and medical scientists to study and simulate ion channel proteins. 

\section*{Acknowledgements}
 This work was partially supported by the National Science Foundation, USA, through award number  DMS-2153376, and the Simons Foundation, USA,  through research award 711776. 
 


\end{document}